\documentclass[11pt, reqno, english]{amsart}  
\usepackage[utf8]{inputenc}
\usepackage[T1]{fontenc}
\usepackage{amsmath,amsthm}
\usepackage{amsfonts,amssymb}
\usepackage{url}
\usepackage{mathtools}  
\usepackage[colorlinks=true,urlcolor=blue,linkcolor=red,citecolor=magenta]{hyperref}
\usepackage{enumerate, paralist}
\usepackage{tikz-cd}
\usepackage[left=1in,right=1in,top=1in,bottom=1in]{geometry}
\usepackage{tikz}
\usepackage{xcolor}
\usetikzlibrary{shapes.geometric, calc}
\definecolor{myblue}{HTML}{0065d6}
\definecolor{mybluetoo}{HTML}{94c4f3}
\definecolor{mygray}{HTML}{d6d6d6}
\usepackage{hyperref}
\usepackage{soul}
\usepackage{mathrsfs}
\usepackage{caption}
\usepackage{fullpage}
\usepackage{placeins}
\usepackage{wasysym}
\usepackage{svg}
\usepackage{import}

\newtheoremstyle{myremark} % name
{7pt}                    % Space above
{7pt}                    % Space below
{}  	                 % Body font
{}                           % Indent amount
{\bf}       	         % Theorem head fuont
{.}                          % Punctuation after theorem head
{.5em}                       % Space after theorem head
{}  % Theorem head spec (can be left empty, meaning ?normal?)

\theoremstyle{plain}
\newtheorem{lemma}{Lemma}[section]
\newtheorem{theorem}[lemma]{Theorem}
\newtheorem{corollary}[lemma]{Corollary}
\newtheorem{proposition}[lemma]{Proposition}

\theoremstyle{definition}
\newtheorem{definition}[lemma]{Definition}
\newtheorem{question}[lemma]{Question}

\theoremstyle{myremark}
\newtheorem{remark}[lemma]{Remark}
\newtheorem{example}[lemma]{Example}

\newcommand{\R}{\mathbb{R}}
\newcommand{\RP}{\mathbb{RP}}
\newcommand{\Z}{\mathbb{Z}}

\newcommand{\cP}{\mathcal{P}}
\newcommand{\cU}{\mathcal{U}}
\newcommand{\cV}{\mathcal{V}}
\newcommand{\cW}{\mathcal{W}}

\newcommand{\diam}{\mathrm{diam}}

\newcommand{\cost}{\mathrm{cost}}

\newcommand{\spread}{\mathrm{spread}}
\newcommand{\pack}{\mathrm{pack}}
\newcommand{\hdim}{\mathrm{hdim}}
\newcommand{\conn}{\mathrm{conn}}
\newcommand{\cat}{\mathrm{Cat}}
\newcommand{\mass}{\mathrm{mass}}
\newcommand{\gr}{\mathrm{Gr}}

\newcommand{\avr}[2]{\mathrm{AVR}(#1;#2)}
\newcommand{\avrm}[2]{\mathrm{AVR}^m(#1;#2)}
\newcommand{\tavr}[2]{\mathrm{TAVR}(#1;#2)}
\newcommand{\tavrm}[2]{\mathrm{TAVR}^m(#1;#2)}

\newcommand{\bor}[2]{\mathrm{Bor}(#1;#2)}

\newcommand{\vr}[2]{\mathrm{VR}(#1;#2)}
\newcommand{\vrm}[2]{\mathrm{VR}^m(#1;#2)}

\newcommand{\im}{\mathrm{im}}

\newcommand{\supp}{\mathrm{supp}}

\newcommand{\norm}[1]{\lVert#1\rVert}

\begin{document}

\title{
Anti-Vietoris--Rips metric thickenings and Borsuk graphs
}

\author{Henry Adams}
\email{henry.adams@ufl.edu}
\author{Alex Elchesen}
\email{alex.elchesen@colostate.edu}
\author{Sucharita Mallick}
\email{sucharitamallick@ufl.edu}
\author{Michael Moy}
\email{michael.moy@colostate.edu}

\subjclass{Primary 55N31,
%Algebraic topology ? Homology and cohomology theories in algebraic topology ? Persistent homology and applications, topological data analysis
Secondary
05C15,
%Combinatorics ? Graph theory ? Coloring of graphs and hypergraphs
51F99,
%: Geometry ? Metric Geometry ? None of the above, but in this section
55R10.
%Algebraic topology ? Fiber spaces and bundles in algebraic topology ? Fiber bundles in algebraic topology
}

\maketitle

\begin{abstract}
For $X$ a metric space and $r\ge 0$, the \emph{anti-Vietoris--Rips metric thickening $\avrm{X}{r}$} is the space of all finitely supported probability measures on $X$ whose support has spread at least $r$, equipped with an optimal transport topology.
We study the anti-Vietoris--Rips metric thickenings of spheres.
We have a homeomorphism $\avrm{S^n}{r} \cong S^n$ for $r > \pi$, a homotopy equivalence $\avrm{S^n}{r} \simeq \RP^{n}$ for $\frac{2\pi}{3} < r \le \pi$, and contractibility $\avrm{S^n}{r} \simeq *$ for $r=0$.
For an $n$-dimensional compact Riemannian manifold $M$, we show that the covering dimension of $\avrm{M}{r}$ is at most $(n+1)p-1$, where $p$ is the packing number of $M$ at scale $r$.
Hence the $k$-dimensional \v{C}ech cohomology of $\avrm{M}{r}$ vanishes in all dimensions $k\geq (n+1)p$.
We prove more about the topology of $\avrm{S^n}{\frac{2\pi}{3}}$, which has vanishing cohomology in dimensions $2n+2$ and higher.
We explore connections to chromatic numbers of Borsuk graphs, and in particular we prove that for $k>n$, no graph homomorphism $\bor{S^k}{r} \to \bor{S^n}{\alpha}$ exists when $\alpha > \frac{2\pi}{3}$.
\end{abstract}

%For $X$ a metric space and $r\ge 0$, the \emph{anti-Vietoris--Rips metric thickening $\mathrm{AVR^m}(X;r)$} is the space of all finitely supported probability measures on $X$ whose support has spread at least $r$, equipped with an optimal transport topology.
%We study the anti-Vietoris--Rips metric thickenings of spheres.
%We have a homeomorphism $\mathrm{AVR^m}(S^n;r) \cong S^n$ for $r > \pi$, a homotopy equivalence $\mathrm{AVR^m}(S^n;r) \simeq \mathbb{RP}^{n}$ for $\frac{2\pi}{3} < r \le \pi$, and contractibility $\mathrm{AVR^m}(S^n;r) \simeq *$ for $r=0$.
%For an $n$-dimensional compact Riemannian manifold $M$, we show that the covering dimension of $\mathrm{AVR^m}(M;r)$ is at most $(n+1)p-1$, where $p$ is the packing number of $M$ at scale $r$.
%Hence the $k$-dimensional \v{C}ech cohomology of $\mathrm{AVR^m}(M;r)$ vanishes in all dimensions $k\geq (n+1)p$.
%We prove more about the topology of $\mathrm{AVR^m}(S^n;\frac{2\pi}{3})$, which has vanishing cohomology in dimensions $2n+2$ and higher.
%We explore connections to chromatic numbers of Borsuk graphs, and in particular we prove that for $k>n$, no graph homomorphism $\mathrm{Bor}(S^k;r) \to \mathrm{Bor}(S^n;\alpha)$ exists when $\alpha > \frac{2\pi}{3}$.

\setcounter{tocdepth}{1}
\tableofcontents

\section{Introduction}

Given a metric space $(X,d)$ and a scale parameter $r>0$, the Vietoris--Rips simplicial complex $\vr{X}{r}$ has vertex set $X$, and as simplices all finite subsets of $X$ of diameter at most $r$.
In applied topology, Vietoris--Rips complexes have been used to find the ``shape'' of datasets using persistent homology~\cite{Carlsson2009}.
In these settings, the dataset is often thought of as a finite sample drawn from a manifold.
Latschev's result~\cite{Latschev2001} supports the use of Vietoris--Rips complexes in this way.
Indeed, Latschev proves that if a dataset $X$ is Gromov--Hausdorff close to a manifold $M$, then $\vr{X}{r}$ recovers the homotopy type of $M$ 
for scale parameters $r$ that are sufficiently small compared to the curvature of the manifold and sufficiently large compared to the density of the sample.

If a dataset $X$ is sampled from a manifold $M$, then as more and more points are drawn, the persistent homology of the dataset $X$ converges to the persistent homology of the manifold $M$.
%What is the persistent homology of a manifold?
Hausmann's result guarantees that at small scales, the Vietoris--Rips complex $\vr{M}{r}$ is homotopy equivalent to the manifold $M$.
But, in persistent homology, the idea is to let the scale grow from small to large, instead of trying to guess which scales are sufficiently small.
What is the persistent homology of a manifold at larger scales?
The only known complete answer is when the manifold $M$ is the circle $S^1$, in which case $\vr{S^1}{r}$ obtains the homotopy types of the circle $S^1$, the 3-sphere $S^3$, the 5-sphere $S^5$, the 7-sphere $S^7$, \ldots as the scale increases~\cite{AA-VRS1}.
This remains a difficult problem for other manifolds, including Vietoris--Rips complexes of higher-dimensional spheres.
The Vietoris--Rips metric thickening $\vrm{X}{r}$, which consists of a finitely-supported probability measures on $X$ whose support has diameter at most $r$, was introduced in~\cite{AAF} in order to advance our understanding of Vietoris--Rips complexes of manifolds.

In a dissimilarity matrix (for example a distance matrix), large values correspond to dissimilar data points, whereas in a similarity matrix, large values correspond to similar data points.
Similarity matrices arise frequently in data science, for example in neuroscience where two neurons with similar spike train data are assigned a large measure of similarity, in time series analyses when similar time series are assigned a high correlation value, and in network science when similar nodes are assigned a large affinity value.
In his PhD thesis~\cite{JeromeRoehm}, Roehm makes the case well that (1) when studying neuron firings, similarity matrices are more natural than dissimilarity matrices, and (2) when data is presented in a similarity matrix, it is more natural to connect points that corresponds to ``higher values'' first, instead of transforming the matrix to a dissimilarity matrix and connecting the points that correspond to lower values first.
This corresponds to building an \emph{anti-Vietoris--Rips simplicial complex}, adding edges with values bounded from below and taking a clique complex (see~\cite{engstrom2009complexes,Jefferson-AATRNtalk2021,JeffersonAntihomology}), instead of a Vietoris--Rips construction with edge values bounded from above.

We initiate the study of anti-Vietoris--Rips (or anti-VR) thickenings of a metric space.
If $X$ is a metric space and $r\ge 0$ is a scale parameter, then the \emph{anti-Vietoris--Rips metric thickening} $\avrm{X}{r}$ consists of all finitely-supported probability measures on $X$ whose support has spread at least $r$.
A subset $\sigma \subseteq X$ has \emph{spread} at least $r$ if $d(x,x') \ge r$ for all $x,x'\in \sigma$ with $x\neq x'$.
If $r$ is larger than the diameter of $X$, then $\avrm{X}{r}$ consists of all Dirac delta measures, and hence is isometric to $X$.
However, as $r$ decreases below the diameter of $X$, the homotopy type of $\avrm{X}{r}$ may change.
We note that the $r$-packing number of $X$ places an upper bound on the number of points that may be in the support of a probability measure in $\avrm{X}{r}$.
Once $r$ has decreased down to zero, then $\avr{X}{0}$ is the space of all finitely supported probability measures in $X$, which is contractible.

We consider in particular the case of anti-VR thickenings of $n$-spheres.
We equip the $n$-sphere $S^n$ with the geodesic metric, so that $S^n$ has diameter $\pi$ 
(but all of our results have analogues if the sphere is instead equipped with the Euclidean metric).
If $r>\pi$, then $\avrm{S^n}{r}$ is isometric to $S^n$.
We prove in Theorem~\ref{thm:avrmSn-homotopy-type} that if $\frac{2\pi}{3} < r\ \le \pi$, then $\avrm{S^n}{r}$ is homotopy equivalent to $\RP^n$.
In Theorem~\ref{thm:cov-dim} we show that for an $n$-dimensional Riemannian manifold $M$, the covering dimension of $\avrm{M}{r}$ is at most $(n+1)p-1$, where the packing number $p=\pack_M(r)$ is the largest number of points that can be placed in $M$ so that any two points are at distance at least $r$ apart.
Thus $\avrm{M}{r}$ has no cohomology in dimensions $(n+1)p$ and higher.
In Theorem~\ref{thm:hom_type_avrm(S^n;2pi/3)}, we also give initial information on the topology of $\avrm{S^n}{\frac{2\pi}{3}}$, which has vanishing cohomology in dimensions $2n+2$ and higher.
The homotopy types of $\avrm{S^n}{r}$ for $0<r<\frac{2\pi}{3}$ remain mysterious.

We will also see that the homotopy types of anti-VR thickenings can be used as obstructions to the existence of certain graph homomorphisms between Borsuk graphs.
Using the homotopy type of $\avrm{S^n}{r}$, we show in Theorem~\ref{thm:no-graph-homomorphism} that graph homomorphisms from a Borsuk graph on a higher-dimensional sphere to a Borsuk graph on lower dimensional sphere cannot exist unless the scale of the codomain is sufficiently relaxed.
Indeed, for all $k>n$ and $r<\pi$, we show there is no graph homomorphism $f\colon \bor{S^k}{r}\to \bor{S^n}{\alpha}$ when $\alpha \geq \frac{2\pi}{3}$.
This improves upon the result implied by what is known about the chromatic numbers of these graphs.

Our paper is organized as follows.
We survey related work in Section~\ref{sec:related-work}, and cover background material in Section~\ref{sec:background}.
In Section~\ref{sec:anti-vr}, we define anti-VR simplicial complexes and metric thickenings.
In Section~\ref{sec:avr-large}, we prove that the anti-VR thickening of the $n$-sphere $S^n$, for scales $\frac{2\pi}{3} < r \leq \pi$, is homotopy equivalent to the $n$-dimensional real projective space $\RP^n$.
In Section~\ref{sec:cov-dim} we bound the covering dimension of anti-VR thickenings of compact manifolds, and hence the highest dimension in which cohomology can appear for these spaces.
In Section~\ref{sec:avr-small} we give some information about the topology of $\avr{S^n}{r}$ at the first scale small enough so that the homotopy type is no longer $\RP^n$.
In Section~\ref{sec:TAVR}, we introduce total anti-VR complexes and thickenings.
We describe a connection to chromatic numbers, especially of Borsuk graphs, in Section~\ref{sec:choromatic-numbers-borsuk-grp}.
Then, using the homotopy equivalences from Sections~\ref{sec:avr-large} and~\ref{sec:TAVR}, we prove a no-graph-homomorphism result for Borsuk graphs in Section~\ref{sec:no-graph-homomorphisms}.
We conclude with a list of open questions in Section~\ref{sec:conclusion}.

\section{Related work}
\label{sec:related-work}

We give a short survey of Vietoris--Rips complexes and metric thickenings, which are well-studied.
After that, we move on to anti-Vietoris--Rips complexes.
We end by surveying some papers that make use of anti-Vietoris--Rips perspectives, even if they may not be directly related.

\subsection*{Vietoris--Rips complexes and thickenings}

Vietoris--Rips complexes are commonly used tools in applied topology for estimating the ``shape'' of data~\cite{Carlsson2009}.
Given a metric space $X$ and a scale $r\ge 0$, the Vietoris--Rips complex $\vr{X}{r}$ contains a simplex for each finite set of $X$ of diameter at most $r$.
They were introduced by Vietoris~\cite{VietorisVRcplx} in order to define a cohomology theory for metric spaces, and used by Rips in geometric group theory (for example, Rips proved that the Vietoris--Rips complex of a $\delta$-hyperbolic group is contractible if the scale is chosen to be at least $4\delta$~\cite{GromovHyperbolicGroups}).

Hausmann proves that for a closed Riemannian manifold $M$, when the scale $r$ is sufficiently small compared to the injectivity radius of $M$, the Vietoris--Rips complex $\vr{M}{r}$ is homotopy equivalent to $M$~\cite{Hausmann}.
Latschev further proves that if a space $X$ is Gromov-Hausdorff close to $M$ and the scale parameter $r$ is in the right range, then $\vr{X}{r}$ recovers the homotopy type of $M$~\cite{Latschev2001}.
In practice, to recover the homology groups of $M$, one looks at the persistent homology of $\vr{X}{r}$ as the scale $r$ varies.
The stability of persistent homology~\cite{PersistenceStability,chazal2009gromov} implies that if $X$ is close to $M$ in the Gromov--Hausdorff distance, then the persistent homology of $\vr{X}{-}$ is close to the persistent homology of $\vr{M}{-}$.

In order to use the persistent homology of $\vr{X}{-}$ to recover the homology of $M$, one typically uses the principle that features that persist for a long range of scales $r$ detect the homology or the homotopy type of $M$~\cite{ChazalOudot2008}.
By contrast, short-lived persistent homology features are indicative of local geometry.
Indeed,~\cite{bubenik2020persistent} shows that 
the short-lived features or ``noise'' can be used to detect curvature of disks from which points have been sampled, and~\cite{robins2000computational,macpherson2012measuring,adams2018fractal,schweinhart2019persistent,schweinhart2020fractal} consider definitions of fractal dimension using persistence homology.
Both types of features, long- and short-lived, have found use in machine learning applications~\cite{bubenik2015statistical,PersistenceImages,bendich2016persistent,hiraoka2016hierarchical,nakamura2015persistent,nakamura2015persistent,adams2021topology}.

Vietoris--Rips metric thickenings $\vrm{X}{r}$ were introduced in~\cite{AAF} because they are often better behaved than Vietoris--Rips simplicial complexes.
In particular, the inclusion $X\hookrightarrow \vrm{X}{r}$ is continuous, whereas the inclusion $X\hookrightarrow \vrm{X}{r}$ need not be continuous for $X$ infinite.
The paper~\cite{AAF} gives the first new homotopy type $S^n*\frac{\mathrm{SO}(n+1)}{A_{n+2}}$ appearing for Vietoris--Rips metric thickenings of the $n$-sphere $S^n$ for all $n$.
The paper~\cite{moy2023vietoris} provides the homotopy types of $\vrm{S^1}{r}$ for all values of the scale parameter $r$, which are again $S^1$, $S^3$, $S^5$, $S^7$,\ldots as the scale parameter increases, now proven using more geometric and Morse-theoretic techniques; see also~\cite{PersistentEquivariantCohomology}.
Furthermore, the algebraic invariants of Vietoris--Rips metric thickenings have been closely related to those of Vietoris--Rips simplicial complexes; see~\cite{MoyMasters,AMMW,HA-FF-ZV,gillespie2022homological,gillespie2024vietoris}.
For a totally bounded metric space $X$, the persistent homology of Vietoris--Rips metric thickenings agrees with the persistent homology of Vietoris--Rips complexes, so long as we ignore whether endpoints of intervals are open or closed~\cite{AMMW,gillespie2024vietoris}.
Hence the persistent homology of Vietoris--Rips and \v{C}ech metric thickenings are stable under perturbations with respect to the Gromov--Hausdorff distance on the underlying metric space $X$.
%The paper~\cite{AAF} gives the first new homotopy type $S^n*\frac{\mathrm{SO}(n+1)}{A_{n+2}}$ appearing for Vietoris--Rips metric thickenings of the $n$-sphere $S^n$ for all $n$.
%The paper~\cite{moy2023vietoris} provides the homotopy types of $\vrm{S^1}{r}$ for all values of the scale parameter $r$, which are $S^1$, $S^3$, $S^5$, $S^7$, \ldots as the scale parameter increases.
%Indeed, the topology of Vietoris--Rips metric thickenings is often better behaved, in part because the inclusion $M\hookrightarrow \vrm{M}{r}$ into the metric thickening is continuous, whereas the inclusion $M\hookrightarrow \vr{M}{r}$ into the simplicial complex is discontinuous for $M$ a manifold of dimension at least one.
%In a similar fashion, in this paper we will study the homotopy types of the \emph{anti}-Vietoris--Rips thickenings on $n$-spheres.

\subsection*{Anti Vietoris--Rips complexes}

The perspective taken with \emph{anti}-Vietoris--Rips simplicial complexes $\avr{X}{r}$ is instead to connect vertices in $X$ when they are ``far apart''.
Indeed, the simplices in $\avr{X}{r}$ are all finite subsets of $X$ that have spread at least $r$, i.e.\ all finite $\{x_0,\ldots,x_k\}\subseteq X$ with $d(x_i,x_j)\ge r$ for all $i\neq j$.
Compared to Vietoris--Rips complexes, there are relatively fewer references on anti-Vietoris--Rips complexes.
We survey the ones that we know about.

In~\cite{engstrom2009complexes}, Engstr{\"o}m proves two main results about anti-VR complexes.
First, in~\cite[Proposition~4.2]{engstrom2009complexes} Engstr{\"o}m gives a recursive formula for the homotopy type of the anti-VR complex of any finite subset of the real line, in terms of wedge sums of suspensions of smaller anti-VR complexes.
Second, in~\cite[Proposition~4.3]{engstrom2009complexes} Engstr{\"o}m proves that if $X\subseteq \Z^2$ has $n$ vertices, then $\avr{X}{r}$ is $\lfloor \frac{n-9}{6} \rfloor$-connected for $1<r<\sqrt{2}$.

Consider an unweighted graph $G$ with the shortest path metric.
The $r$-th graph power $G^{(r)}$ is a graph which has vertices $V(G)$ and two vertices are adjacent if the distance between them is at most $r$ in the shortest path metric~\cite{engstrom2009complexes}.
Note that $G^{(1)}=G$.
Also, notice that the independence complex of the $r$-th graph power $I(G^{(r)})$ is the same as $\avr{G}{r+1}$ (and so in particular, $I(G)=\avr{G}{2}$).
For further references on independence complexes, see for example~\cite{berghoff2020homology,ehrenborg2006topology,engstrom2008independence,Kozlov}, including filtered variants~\cite{abdelmalek2023chordal,deshpande2022distance}.

In statistical physics, solutions to certain hard-particle models can be translated to the enumeration of independent sets of a certain square grid graph $G$.
Looking at the hard-particle model at ``activity $-1$'' on these graphs is the same as the difference of the number of odd and even cardinality independent sets.
This is the same as computing the Euler characteristic of $\avr{G}{2}$; see~\cite{IndCplxSquareGrids} for further details.

We refer the reader to~\cite{Jefferson-AATRNtalk2021} for a video by Jefferson introducing anti-VR complexes; see also~\cite{JeffersonAntihomology}.
From a dataset (a finite sample of points from a circle, an ellipse, a square, etc.), Jefferson forms an anti-Vietoris--Rips complex.
Then, the author computes the persistent homology of filtrations given by these anti-VR complexes as the scale is relaxed.

In his PhD thesis~\cite{JeromeRoehm}, Roehm studies anti-geometric persistence.
He discusses dissimilarity matrices and computes barcodes from them, which he calls ``anti-barcodes''.
An example of a dissimilarity matrix would be any distance matrix --- large entries correspond to dissimilar data points.
But, in a dissimilarity matrix, you do not necessarily assume that the triangle inequality is satisfied.
By contrast, in a \emph{similarity} matrix, large entries correspond to similar data points (think highly correlated data points).
For example, when measuring neuron firings, it might be most natural to first construct a notion of similarity which is the correlation between two firing patterns.
Roehm makes the point well that when data is presented in the form of a similarity matrix, you might want to proceed directly and construct the anti-VR complex (which will connect similar data points first), instead of first transforming your similarity matrix into a dissimilarity matrix\footnote{
There are several ways (such as $x\mapsto \frac{1}{x}$ or $x\mapsto e^{-x}$) but no canonical way to turn a similiarity matrix into a dissimilarity matrix.
Also, one can turn a dissimilarity matrix into a nearby metric space using multidimensional scaling (MDS) or related techniques; see~\cite{bibby1979multivariate,cox2000multidimensional,groenen2014past,MDScircle,de2022component}.
}
and then computing a Vietoris--Rips complex.
In \cite[Theorem~4.5.1]{JeromeRoehm}, Roehm shows given any two arbitrary persistence barcodes $B_1$ and $B_2$, there exists a dissimilarity matrix whose anti-Vietoris--Rips complexes yield barcode $B_1$ and whose Vietoris--Rips complexes of the transposed similarity matrix yield barcode $B_2$; hence in general there is no duality between what is recovered by the anti-VR and the VR constructions.
Roehm studies anti-Vietoris--Rips complexes of certain finite metric spaces, such as evenly spaced points on a line, evenly spaced points on a circle, or a finite subset of a sphere.

By contrast, we focus our attention on anti-Vietoris--Rips metric thickenings of \emph{infinte} metric spaces, such as the entire $n$-sphere.
Though Vietoris--Rips metric thickenings of infinite metric spaces (such as spheres) are closely related to Vietoris--Rips simplicial complexes of finite samples thereof, it is our impression that anti-Vietoris--Rips metric thickenings of infinite metric spaces are more related to \emph{total} anti-Vietoris--Rips simplicial complexes of finite samples thereof; see Section~\ref{sec:TAVR}.

We note that if $X$ is an unbounded metric space, then $\avr{X}{r}$ is contractible for any $r>0$.
Indeed, we show that $\pi_i(\avr{X}{r})$ is trivial for all $i\ge 0$.
By compactness, any map $f\colon S^i \to \avr{X}{r}$ lands in a finite subcomplex with finite vertex set $V$.
Since $V$ is finite, there is some vertex $x\in X\setminus V$ at distance greater than $r$ from each vertex of $V$.
Therefore, the image of $f$ lives in a cone on vertex $x$, and hence $f$ is nullhomotopic.
By Whitehead's theorem, it follows that $\avr{X}{r}$ is contractible.
So, when studying anti-Vietoris--Rips complexes, it is most interesting to restrict attention to bounded metric spaces.

Chapter~5 of the fourth author's PhD thesis~\cite{MoyThesis} is a follow-up project based on this paper, even though it appeared earlier.
In the case of $S^1$, Moy proves that the homotopy types of $\avrm{S^1}{r}$ are $S^1$, $S^1$, $S^3$, $S^3$, $S^5$, $S^5$, $S^7$, $S^7$,\ldots as one decreases the scale parameter, until $\avrm{S^1}{r}$ is finally contractible when $r=0$.
Why is each odd-dimensional sphere listed twice?
One has $\avrm{S^1}{r}\simeq S^1$ for $\frac{2\pi}{3}<r$, although for $\frac{2\pi}{3} < r' \le \pi < r$ the inclusion $\avrm{S^1}{r}\hookrightarrow\avrm{S^1}{r'}$ yields a degree two map $S^1\to S^1$.
This degree two map can be seen as the natural 2-fold covering map $S^1 \to \RP^1 = S^1$ (compare Theorem~\ref{thm:avrmSn-homotopy-type}).
More generally for $k \geq 2$, $\avrm{S^1}{r}\simeq S^{2k-1}$ for $\frac{2\pi}{2k+1} < r \le \frac{2\pi}{2k-1}$, and when $\frac{2\pi}{2k+1} < r' \le \frac{2\pi}{2k} < r \le \frac{2\pi}{2k-1}$, the inclusion $\avrm{S^1}{r}\hookrightarrow\avrm{S^1}{r'}$ yields a degree two map $S^{2k-1}\to S^{2k-1}$.

\subsection*{Other related papers}

In~\cite{sankar2022homotopy} Sankar uses a notion of discrete homotopy in order to study questions (and threshold parameters) related to the existence of graph homomorphisms.
In one such example, the vertices of a discrete path lie in an $n$-sphere, and adjacent vertices on the path are required to be nearly antipodal, i.e.\ within an $\varepsilon$ error bound of being antipodal.
Note that each step in such a path moves along an edge of the anti-VR complex of the sphere at scale $\pi-\varepsilon$.
For this reason, this discrete homotopy theory feels related to anti-VR type thickenings.
% \note{Cite the Borsuk--Ulam theorem and related generalizations thereof.}
See~\cite{rieser2021cech,rieser2020semiuniform,bubenik2024homotopy,dochtermann2009hom,barcelo2001foundations,babson2004homotopy,BarceloCapraroWhite,carranza2022cubical} for other variants of discrete homotopy theories.

We ask if there is a relationship between anti-Vietoris--Rips complexes and configurations spaces of balls of positive radii; see for example~\cite{carlsson2012computational,baryshnikov2014min,alpert2021configuration}.
Indeed, whenever the open balls $\{B(x_i;\frac{r}{2})\}_{i=0}^k$ are disjoint in a geodesic space $X$, then $\{x_0,\ldots,x_k\}$ will be a simplex in the anti-Vietoris--Rips complex $\avr{X}{r}$.

We also refer the reader to Barmak's paper~\cite{BarmakStar2013}, especially Theorem~6.2 and Corollary~6.8 within.
Theorem~6.2 gives the bound $\chi(G)\ge \cat(I(G))+1$.
Here the \emph{strong Lusternik--Schnirelmann category $\cat(Y)$} of a space $Y$ is one less that the minimum number of contractible subcomplexes which are needed to cover some CW complex homotopy equivalent to $Y$.
It is known that $\cat(Y)$ is bounded from below by the cup-length of $Y$.
Since the independence complex of a graph is the clique complex of the graph complement, 
the independence complex of the Borsuk graph is a Vietoris--Rips complex of a sphere.
This provides additional motivation for understanding the homotopy types and Lusternik--Schnirelmann categories of Vietoris--Rips complexes of spheres.
That said, Corollary~6.8 of~\cite{BarmakStar2013} shows that the clique number $\omega(G)$ is bounded in-between ($\chi(G)\ge \omega(G)\ge \cat(I(G))+1$), meaning we do not expect this bound to be tight in many cases.

\section{Background and notation}
\label{sec:background}

In this section we define existing concepts and set notation for metric spaces, simplicial complexes, optimal transport, metric thickenings, and fiber bundles.

Let $X$ and $Y$ be topological spaces.
A \emph{map} $f\colon X\to Y$ is a continuous function.
Also, let $I=[0,1]$ be the unit interval.

\subsection*{Metric spaces and spread}

Let $X$ be a metric space with distance function $d\colon X\times X\to \R$.
We let $B(x,r)=\{x'\in X~|~d(x,x')<r\}$ denote the ball of radius $r>0$ centered at $x\in X$.

For $A \subseteq X$, we define the \emph{diameter} of $A$ to be $\diam(A)\coloneqq \sup\{d(x,x')~|~x,x'\in A\}$.
The diameter is zero if $A$ is a singleton.

Similarly, we define the $spread$ of $A$ to be
\[\spread(A) \coloneqq \inf\{d(x,x')\mid x\neq x',\ x,x'\in A\}.\]
In the case when $A$ is a singleton set, we are taking the infimum of the empty set, and thus the spread of a singleton set is $\infty$.
In this paper, we will mostly be considering the spread of finite subsets.

A metric space $(X,d)$ is said to be \emph{totally bounded} if for every $\varepsilon>0$, there exist a finite set of points $F=\{x_1,x_2,\cdots,x_n\}\subseteq X$ such that $X\subseteq \bigcup_{i=1}^n B(x_i,\varepsilon)$.
The inclusion $X\subseteq \bigcup_{i=1}^n B(x_i,\varepsilon)$ means that every point $x\in X$ is less than $\varepsilon$ away from $F$, i.e., $\min_i{d(x,x_i)}<\varepsilon$.
Such a set $F$ is called an \emph{$\varepsilon$-net} on $X$.

\subsection*{Simplicial complexes and Vietoris--Rips complexes}

Let $K$ be a simplicial complex on vertex set $V$.
Even if $V$ is infinite, we remind the reader that a simplex $\sigma\in K$ is required to be a finite subset of $V$.
We identify abstract simplicial complexes with their geometric realizations.

Given a metric space $(X,d)$ and $r>0$, the Vietoris--Rips complex $\vr{X}{r}$ is the simplicial complex with vertex set $X$ and with simplices all finite subsets of $X$ of diameter at most $r$.

\subsection*{Optimal transport distances}
We follow~\cite{villani2003topics}.
Let $(X,d)$ be a metric space.
A Borel measure $\mu$ is said to be a \emph{Radon measure} if it is inner regular, i.e., $\mu(B)=\sup\{\mu(K)\mid K \text{ compact, } K\subseteq B\}$, and locally finite, i.e., for every $x\in X$ there exist a neighborhood $U$ of $x$ such that $\mu(U)<\infty$.
Let $\mathcal{P}(X)$ be set of all probability Radon measures such that there exists some $y\in X$ with $\int_X d(x,y)d\mu<\infty$.

The Wasserstein or Kantorovich~\cite{vershik2013long} metric is a metric defined on $\mathcal{P}(X)$.
Given $\mu,\nu\in \mathcal{P}(X)$, let $\Pi(\mu,\nu)\subseteq \mathcal{P}(X\times X)$ be the set of all probability Radon measures $\pi$ with marginals $\mu$ and $\nu$, i.e., for any Borel set $E\subseteq X$ we have $\pi(E\times X)=\mu(E)$ and $\pi(X\times E)=\nu(E)$.
We refer to $\pi\in\Pi(\mu,\nu)$ as a \emph{transport plan} between $\mu$ and $\nu$.
To see that the set $\Pi(\mu,\nu)$ is non-empty, consider the product measure $\mu\times \nu\in \Pi(\mu,\nu)$.
For $\mu,\nu\in \mathcal{P}(X)$, the \emph{$1$-Wasserstein} (or \emph{optimal transport}) metric $W_1$ is defined as
\[W_1(\mu,\nu)\coloneqq\inf_{\pi\in \Pi(\mu,\nu)}\int_{X\times X} d(x,y)d\pi.\]
Here, $\int_{X\times X} d(x,y)d\pi$ is the \emph{cost} associated to the transport plan $\pi$, which we denote by $\cost(\pi)$.
The infimum in this definition is attained and the infimum gives the solution to the Monge--Kantorovich problem.
The $1$-Wasserstein distance is therefore the cost of an optimal transport plan.

In a sense, the optimal transport metric $W_1$ generalizes the metric $d$ in $X$.
For any $x\in X$, let $\delta_x$ be the Dirac delta measure supported on $x$, so $\delta_x\in \mathcal{P}(X)$.
Then $x\mapsto \delta_x$ is an isometry, i.e., for any $x,y \in X$, $W_1(\delta_x,\delta_y)=d(x,y)$.

If $\pi_1$ and $\pi_2$ are two transport plans between $\mu$ and $\nu$, then for any $0\leq \lambda\leq 1$, so is $\lambda\pi_1+(1-\lambda)\pi_2$.
Also, 
\[\cost(\lambda\pi_1+(1-\lambda)\pi_2)=\lambda \cost(\pi_1)+(1-\lambda)\cost(\pi_2).\]

We will find use for partial transport plans in Section~\ref{sec:cov-dim}, indeed to prove Theorem~\ref{thm:cov-dim}.
A \emph{partial transport plan} between $\mu, \nu\in \mathcal{P}(X)$ is a Radon measure $\sigma$ such that for any Borel set $E\subseteq X$, we have $\sigma(E\times X)\leq\mu(E)$ and $\sigma(X\times E)\leq\nu(E)$.
The \emph{mass} of a partial transport plan $\sigma$ is $\mass(\sigma)=\int_{X\times X}d\sigma$.
The mass is always less than or equal to one, and strict inequality is possible.
If a partial transport plan has mass one, then it is a transport plan.
The \emph{cost} of the partial transport plan $\sigma$ is $\int_{X\times X} d(x,y)d\sigma$.
Though partial transport plans are relevant for the optimal partial transport problem (see~\cite{caffarelli2010free,figalli2010optimal}), we will not consider that problem here.
Instead, we use partial transport plans to bound the cost of (full) transport plans in spaces of bounded diameter, as per the following lemma.

\begin{lemma}
\label{lem:ExtnOfOptimalTransport}
Let $(X,d)$ be a metric space with $\diam(X)<\infty$, let $\mu,\nu\in \mathcal{P}(X)$, and let $\sigma$ be a partial transport plan between $\mu$ and $\nu$.
Then $\sigma$ can be extended to a (full) transport plan $\tilde{\sigma}$ between $\mu$ and $\nu$ with
\[\cost(\tilde{\sigma})\leq \cost(\sigma)+ \left( 1-\mass(\sigma)\right)\diam(X).\]
\end{lemma}

See Appendix~\ref{app:leftover-proofs} for the proof.

\subsection*{Metric thickenings and Vietoris--Rips metric thickenings}
Given a metric space $(X,d)$ and $r \ge 0$, the \emph{Vietoris--Rips metric thickening $\vrm{X}{r}$} is the space consisting of probability measures $\sum_{i=1}^k\lambda_i\delta_{x_i}$ with $\lambda_i > 0$, with $\sum_{i=1}^k\lambda_i=1$, with $x_i\in X$ for all $i$, and with $\diam\{x_1, \dots, x_k\} \leq r$~\cite{AAF}.
Here, $\delta_{x_i}$ is a Dirac delta measure at $x_i$ and the set $\{x_1,\ldots, x_k\}$ is said to be the \emph{support} of $\sum_{i=1}^k\lambda_i\delta_{x_i}$.
The space $\vrm{X}{r}$ is equipped with the optimal transport topology.
For convenience, we use the $1$-Wasserstein distance $W_1$, although this choice is not critical since the $p$-Wasserstein distance $W_p$ induces the same topology for all $1\le p<\infty$~\cite{bogachev2018weak,gibbs2002choosing}.

The homotopy types of Vietoris--Rips simplicial complexes and Vietoris--Rips metric thickenings are tightly connected.
If $X$ is finite then $\vr{X}{r}$ and $\vrm{X}{r}$ are homeomorphic (\cite[Proposition~6.2]{AAF}; see also~\cite[Equation~1.B(c)]{Gromov}).
Adams, Moy, M\'{e}moli, and Wang proved that for $X$ totally bounded, the persistent homology modules of $\vr{X}{r}$ and $\vrm{X}{r}$ are $\varepsilon$-interleaved for any $\varepsilon>0$~\cite{AMMW,MoyMasters}.
Next, Adams, Virk, and Frick proved that if $X$ is totally bounded, then $\vr{X}{r}$ and $\vrm{X}{r}$ (both with the $<$ convention) have isomorphic homotopy groups~\cite{HA-FF-ZV}, but they do not show that they are weakly homotopy equivalent (they do not find a single map inducing these isomorphisms).
Gillespie proved that they have isomorphic homology groups~\cite{gillespie2022homological}.
This sequence of results was notably strengthened in~\cite{gillespie2024vietoris}, when Gillespie proved that the natural continuous map $|\vr{X}{r}|\to\vrm{X}{r}$ (both with the $<$ convention) defined by $\sum_i \lambda_i x_i \mapsto \sum_i \lambda_i \delta_{x_i}$ induces an isomorphism on all homotopy groups, and hence is a weak homotopy equivalence.
This implies that the persistent homology modules of $\vr{X}{r}$ and $\vrm{X}{r}$ are not only $\varepsilon$-interleaved for any $\varepsilon>0$, but are furthermore isomorphic.

\subsection*{Fiber bundles and fibrations}
A \emph{fiber bundle} is a tuple $(E,B,\pi,F)$, where $E$, $B$, and $F$ are topological spaces and $\pi\colon E\to B$ is a continuous map satisfying the following \emph{local triviality condition}: 
for each $b\in B$, there exists a neighborhood $U$ of $b$, and a homeomorphism $h\colon \pi^{-1}(U)\to U\times F$  such that the following diagram commutes, where $p_1$ is the projection to the first coordinate:
\begin{center}
    \begin{tikzcd}[row sep=large, column sep=large]
   \pi^{-1}(U) \ar[d,swap,"\pi"] \ar[r,"h"] & U\times F \ar[dl,"p_1"] \\
    U 
  \end{tikzcd}
\end{center}
Here $E,B$, and $F$ are called the \emph{total space}, \emph{base space}, and \emph{fiber}, respectively.

A mapping $\pi:E\to B$ between topological spaces $E$ and $B$ satisfies the \emph{homotopy lifting property} for a space $X$ if, for every homotopy $H\colon X\times I\to B$ and for every lift $\tilde{H}_0$ of $H_0 \coloneqq H|_{X\times 0}$ to $E$ (i.e., $\pi \circ \tilde{H}_0 = H_0$), there exists a lift $\tilde{H}\colon X\times I\to E$ of $H$ (i.e., $\pi\circ \tilde{H} = H$) which satisfies $\tilde{H}|_{X\times 0} = \tilde{H}_0$.

\begin{equation} \label{eq:homotopy_lifting_property}
  \begin{tikzcd}
    X\times \{0\} \ar[r,"\tilde{H}_0"] \ar[d,swap,hook,""]& E \ar[d,"\pi"] \\
    X\times I \ar[r,"H"] \ar[ur,dashed,"\tilde{H}"] & B
  \end{tikzcd}
\end{equation}

A mapping $\pi:E\to B$ is called a \emph{(Hurewicz) fibration} if it satisfies the homotopy lifting property for every space $X$.
It is called a \emph{Serre fibration} if it satisfies the homotopy lifting property for all CW complexes.

\section{Anti-VR complexes and metric thickenings}
\label{sec:anti-vr}

The main goal of this section is to introduce anti-Vietoris--Rips metric thickenings.
We first define anti-Vietoris--Rips simplicial complexes, then introduce anti-Vietoris--Rips metric thickenings, and then explain how anti-VR thickenings are related to more general metric thickenings.

\subsection{Anti-VR complexes}
We begin with the previously-introduced anti-Vietoris--Rips simplicial complexes.

\begin{definition}
\label{def:AVR}
Let $(X,d)$ be a metric space and $r\ge 0$.
The \emph{anti-Vietoris-Rips (or anti-VR) complex} is the simplicial complex $\avr{X}{r}$ with
\begin{itemize}
\item vertex set $X$
\item simplices of the form $\{x_0,...,x_k\}$ such that $d(x_i,x_j) \ge r$ for all $0 \le i,j \le k$ with $i \neq j$.
\end{itemize}
\end{definition}

Though other authors have used the $>r$ convention (see for example~\cite{engstrom2009complexes}), we use the $\ge r$ convention in this paper.

\subsection{Anti-VR metric thickenings}
\label{ssec:anti-vrm}

We let $\avrm{X}{r}$ denote the associated \emph{anti-VR metric thickening}:

\begin{definition}
\label{def:avrm}
Let $(X,d)$ be a metric space and $r\geq 0$.
The \emph{anti-Vietoris--Rips (or anti-VR) metric thickening $\avrm{X}{r}$} is the set of all finitely-supported probability measures on $X$ whose support has spread at least $r$.
That is,
\begin{align*}
\avrm{X}{r}&=\Big\{\sum_{i=0}^k\lambda_i\delta_{x_i}\mid k\geq 0,\,  \lambda_i> 0,\,\sum_i\lambda_i=1, \,d(x_i,x_j)\geq r\ \forall i\neq j\Big\} \\
&=\Big\{\sum_{i=0}^k\lambda_i\delta_{x_i}\mid k\geq 0,\,  \lambda_i> 0,\,\sum_i\lambda_i=1,\, \spread(\{x_0,\ldots,x_k\})\geq r \Big\},
\end{align*}
equipped with the 1-Wasserstein distance.
\end{definition}

When $r > \diam(X)$, then the only measures whose support have spread at least $r$ are the Dirac delta measures.
Therefore, when $r > \diam(X)$, the anti-VR metric thickening $\avrm{X}{r}$ is the set of all Dirac delta measures in $X$, and hence is isometric to $X$.
When $r=0$, the anti-VR metric thickening $\avrm{X}{r}$ is the space of all finitely supported measures in $X$, which is contractible (via a linear homotopy to any point in this space).
But, for $0 < r < \diam(X)$, the homotopy types of $\avrm{X}{r}$ are in general not well-understood.

Let $(X,d)$ be a metric space.
A subset $S\subseteq X$ is said to be an \emph{$r$-packing} of $X$ if $d(x,x') \ge r$ for all $x,x'\in S$ with $x\neq x'$.
The \emph{packing number}, $\pack_X(r)$, is the largest cardinality of an $r$-packing $S$ of $X$.
The packing number $\pack_X(r)$ of a metric space $X$ determines the maximum number of points in the support of a measure in $\avrm{X}{r}$.

As explained in Section~\ref{sec:background}, the topology of Vietoris--Rips simplicial complexes and metric thickenings are closely related.
In the case when the metric space $X$ is finite, the anti-VR complex $\avr{X}{r}$ and the anti-VR metric thickening $\avrm{X}{r}$ are still homeomorphic by~\cite[Proposition~6.2]{AAF}.
However, when $X$ is infinite, we believe that the topologies of $\avr{X}{r}$ and of $\avrm{X}{r}$ will often diverge considerably.
Indeed, suppose that $X$ is a Riemannian manifold with a finite packing number $p\coloneqq\pack_X(r)$.
Then $\avr{X}{r}$ will be a simplicial complex of dimension $p-1$.
It will typically be the case that $\avr{X}{r}$ has uncountably many (maximal) $(p-1)$-dimensional simplices.
In this setting, we expect the rank of the $(p-1)$-dimensional homology $H_{p-1}(\avr{X}{r})$ of the anti-VR simplicial complex to be much larger than the more natural homology $H_{p-1}(\avrm{X}{r})$ of the anti-VR metric thickening.
As an example, the simplicial complex $\avr{S^1}{\frac{2\pi}{3}}$ has an uncountable number of maximal $2$-dimensional simplices.
We expect that $H_2(\avr{S^1}{\frac{2\pi}{3}})$ has uncountable rank (though this requires further exploration), even though we show in Theorem~\ref{thm:avrmSn-homotopy-type} that $\avrm{S^1}{\frac{2\pi}{3}}\simeq S^3$ and hence $H_2(\avr{S^1}{\frac{2\pi}{3}})=0$ is much smaller.

We next show, in Lemma~\ref{lem:compact}, that the anti-VR metric thickening of a compact metric space is compact.
This is in stark contrast to Vietoris-Rips metric thickenings, which need not preserve compactness (see~\cite[Remark~3.4 and Appendix~D]{AAF}).

\begin{lemma}
\label{lem:compact}
Let $r>0$.
If $(X,d)$ is a compact metric space, then $\avrm{X}{r}$ is compact.
\end{lemma}

\begin{proof}
Since $X$ is compact, $(\mathcal{P}(X),W_1)$ is also compact~\cite[Remark 6.19]{VillaniOldNew}.
%\footnote{\url{https://math.stackexchange.com/questions/4267724/is-the-space-of-probability-measures-on-a-compact-set-is-compact-w-r-t-wasserste}
We will show that $\avrm{X}{r}$ is closed in $\mathcal{P}(X)$.
Let $(\mu_n)$ be a sequence in $\avrm{X}{r}$ converging to some $\mu\in \mathcal{P}(X)$.
We will show that $\mu\in \avrm{X}{r}$, i.e., that $\spread(\supp(\mu))\geq r$.
Suppose to the contrary that there exist distinct points $x,y\in \supp(\mu)$ with $0< d(x,y) < r$.
For any $\varepsilon > 0$, there is an $N$ such that $W_1(\mu, \mu_N) < \varepsilon \min\{ \mu(B_\varepsilon(x)), \mu(B_\varepsilon(y)) \}$, where the right side is necessarily positive because $x,y \in \supp(\mu)$.
Then $\mu_N$ must have a support point $x_N$ within $2 \varepsilon$ of $x$, since otherwise a transport plan between $\mu$ and $\mu_N$ would need to move mass $\mu(B_\varepsilon(x))$ at least distance $\varepsilon$, costing more than $W_1(\mu, \mu_N)$.
Similarly, $\mu_N$ must have a support point $y_N$ within $2 \varepsilon$ of $y$.
Therefore, $d(x,y) \geq d(x_N, y_N) - d(x,x_N) - d(y,y_N) \geq r - 4 \varepsilon$, and since this holds for any $\varepsilon > 0$, we can conclude $d(x,y) \geq r$.
\end{proof}

%\alex{Here's a rigorous proof/explanation if we want to include it.
%It uses KR duality, but we could quickly add that to the introduction and cite Villani:
%By Kantorovich-Rubinstein duality, $W_1(\mu,\mu_N)\geq \int_X fd\mu - \int_X f\mu_N$ for every $1$-Lipschitz function $f:X\to \R$.
%Let $f:X\to \R$ be given by $z\mapsto \varepsilon\max(0,1  - \frac{d(z, B_\varepsilon(x))}{\varepsilon})$\footnote{$f$ is the $1$-Lipschitz function which takes constant value $\varepsilon$ on $B_\varepsilon(x)$, decays with Lipschitz constant $1$ to $0$, and then is constantly $0$ outside of $B_{2\varepsilon}(x)$}.
%Then $f$ is $1$-Lipschitz with support contained in $B_{2\varepsilon}(x)$ and satisfies $\int_X f d\mu \geq \epsilon \mu(B_\varepsilon (x))$.
%Suppose that $\mu_N$ has no support points within $B_{2\varepsilon}(x)$.
%Then $\int_X fd\mu_N = 0$ so that $W_1(\mu,\mu_N) \geq \int_X fd\mu - \int_X f\mu_N = \int_X fd\mu \geq \varepsilon \mu(B_\varepsilon (x))$, contradicting the choice of $N$.
%A similar argument with $y$ replacing $x$ shows that $\mu_N$ contains a support point within $2\varepsilon$ of $y$ as well.}

We will use Lemma~\ref{lem:compact} in Theorem~\ref{thm:cov-dim} when we bound the covering dimension of anti-VR metric thickenings of compact manifolds.

\subsection{Metric thickenings}

More generally, we can define the metric thickening of any simplicial complex whose vertex set is a metric space~\cite{AAF}.

\begin{definition}
\label{def:metric_thickening}
Let $X$ be a metric space and let $K$ be a simplicial complex with vertex set $X$. 
The metric thickening $K^m$ is the following subspace of $\mathcal{P} (X)$, equipped with the restriction of the 1-Wasserstein metric:
\[K^m= \Big\{\sum_{i=0}^k\lambda_i\delta_{x_i}\in \mathcal{P} (X) \mid k\geq 0,\,  \lambda_i> 0,\,\sum_i\lambda_i=1, \{x_0,\ldots,x_k\} \in K \Big\}.\]
\end{definition}

In particular, $\avrm{X}{r}\coloneqq (\avr{X}{r})^m$.

\begin{lemma}
\label{lem:induced-continuous}
Let $X$, $Y$ be compact metric spaces, and let $K$, $L$ be simplicial complexes with vertex sets $V(K)=X$, $V(L)=Y$.
Let $f\colon X\to Y$ be a map of metric spaces such that the induced map $\tilde{f} \colon K^m \to L^m$ on metric thickenings exists.
If $f$ is continuous then~$\tilde{f}$ is continuous.
\end{lemma}

\begin{proof}
%\note{``We prove the case when the metric spaces $X$ and $Y$ are compact (which implies they are Polish); see Appendix~\ref{} for the proof of the general case.''}
Since every compact metric space is a Polish space and in a Polish space $W_1$ convergence is equivalent to weak convergence \cite[Theorem 6.9]{VillaniOldNew}), we do the following.
We need to show that $W_1(\tilde{f}(\mu_n),\tilde{f}(\mu))\to 0$ whenever $W_1(\mu_n,\mu)\to 0$.
For this, it suffices to show that $\int_Y gd\tilde{f}(\mu_n)\to \int_Y gd\tilde{f}(\mu)$ for every bounded continuous function $g:Y\to \R$.
Now since $g\circ f:X\to \R$ is bounded,
\[\int_Y gd\tilde{f}(\mu_n) = \int_X g\circ f d\mu_n \to \int_X g\circ f d\mu = \int_Y gd\tilde{f}(\mu). \qedhere\]
\end{proof}
The compactness assumption on $X$ and $Y$ can be relaxed.
See Lemma~\ref{lem:gen-induced-continuous} in Appendix~\ref{app:gen-induced-continuous} for a more general result.

\begin{corollary}
\label{cor:continuous-on-finite-subset}
Let $X$ and $Y$ be two metric spaces with $X$ finite.
Let $K$ and $L$ be simplicial complexes with vertex sets $V(K)=X$ and $V(L)=Y$.
If $f\colon X\to Y$ is a map such that the induced map $\Tilde{f}\colon K^m\to L^m$ exists, then $\Tilde{f}$ is continuous.
\end{corollary}

\begin{proof}
Since the finite space $X$ is compact and a function from a finite metric space to another space is always continuous, $\im(f)$ is also compact.
Let $L'\subset L$ be the induced subcomplex with vertex set $\im(f)$, meaning $L'=\{\sigma \subseteq \im(f) \mid \sigma \in L\}$.
By Lemma~\ref{lem:induced-continuous}, $\tilde{f}\colon K^m \to (L')^m \hookrightarrow  L^m$ is continuous.
\end{proof}

We use Lemma~\ref{lem:induced-continuous} and Corollary~\ref{cor:continuous-on-finite-subset} in the proof of Theorem~\ref{thm:no-graph-homomorphism}.

The following lemma on continuity of homotopies will be used in the proofs of Lemma~\ref{lem:piCase} and Proposition~\ref{prop:TAVR-to-AVR}.

\begin{lemma}[Lemma~3.9 of~\cite{AAF}]
\label{lem:3.9ofAAR}
Suppose $K^m$ is a metric thickening and $f\colon K^m \to K^m$ is a continuous map such that $H\colon K^m\times [0,1]\to K^m$ given by $H(\mu,t)=(1-t)\mu+tf(\mu)$ is well-defined.
Then $H$ is continuous.
\end{lemma}

\section{Anti-VR thickenings of spheres with large thresholds}
\label{sec:avr-large}

We connect the anti-VR metric thickenings of spheres at large scales to projective spaces, first directly, and afterwards using fiber bundles.

\subsection{Homotopy type of $\avrm{S^n}{r}$ for $\frac{2\pi}{3}< r$}

Let $S^n$ be the $n$-sphere with geodesic metric and with diameter $\pi$.
Let real projective space $\RP^n=S^n/(x\sim -x)$ be the quotient space obtained by identifying antipodal points.
For $x\in S^n$, we let $[x]\in\RP^n$ denot the corresponding equivalence class.
The purpose of this subsection is to prove the following theorem.

\begin{theorem}
\label{thm:avrmSn-homotopy-type}
$\avrm{S^n}{r}\simeq \RP^n$ for all $\frac{2\pi}{3}<r\leq \pi$.
\end{theorem}

\begin{figure}[h]
\def\svgwidth{5.5in}
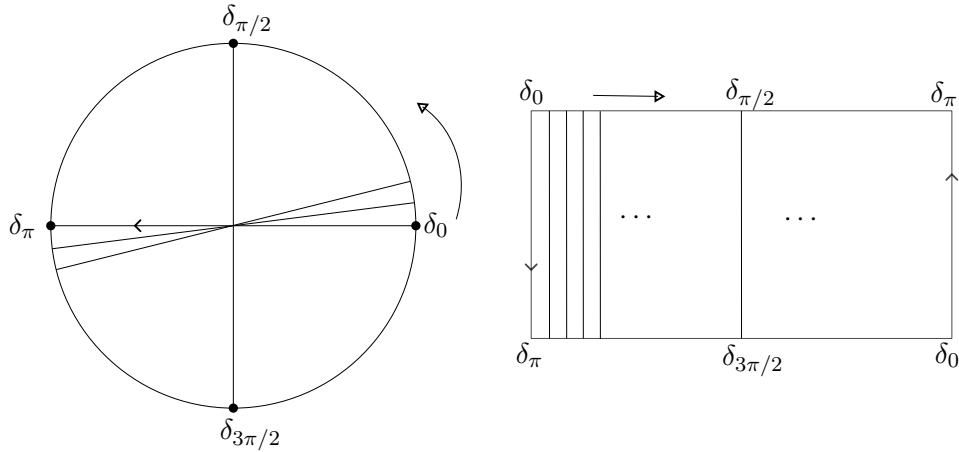
\caption{$\avrm{S^1}{\pi}$ (\emph{left}) is homeomorphic to the M\"{o}bius band (\emph{right}).}
\label{fig:S1-mobius}
\end{figure}

In the case $r=\pi$ and $n=1$, the anti-VR metric thickening $\avrm{S^1}{\pi}$ is homeomorphic to the M\"{o}bius band (see Figure~\ref{fig:S1-mobius}).
The deformation retract we describe in our proof is the well-known deformation retract of the M\"{o}bius band to its central core, which in this setting is the circle $\{\frac{1}{2}\delta_x + \frac{1}{2} \delta_{-x}~|~x\in S^1\}$.

Similarly, $\avrm{S^n}{\pi}$ is homeomorphic to the ``M\"{o}bius $\RP^n$'' (the total space of a nontrivial fiber bundle over $\RP^n$ where each fiber is a unit interval), which deformation retracts to the projective space $\RP^n$.

Before proving Theorem \ref{thm:avrmSn-homotopy-type}, we will need the following results.

\begin{proposition}
\label{CentralCore}
The map $\phi\colon \RP^n \to \avrm{S^n}{\pi}$ defined by $\phi([x]) = \frac{1}{2}\delta_x + \frac{1}{2} \delta_{-x}$ is an embedding of $\RP^n$ into $\avrm{S^n}{\pi}$.
\end{proposition}
\begin{proof}
To see that this is an embedding, consider the map $\phi':S^n\to \avrm{S^n}{\pi}$ given by $\phi'(x) = \frac{1}{2}\delta_x + \frac{1}{2} \delta_{-x}$.
This map satisfies $\phi'(x) = \phi'(-x)$ for all $x\in S^n$ and $\phi\circ p = \phi'$, where $p\colon S^n\to \RP^n$ is the quotient map.
This shows that $\phi$ is continuous.
Since $\phi$ is injective, since $\RP^n$ is compact, and since $\avrm{S^n}{\pi}$ is Hausdorff, we conclude $\phi$ is a homeomorphism onto its image,~\cite[Theorem~26.6]{munkres}.
Thus, $\phi$ is an embedding.
\end{proof}

For $0<r<\pi$, there exists a natural inclusion $\avrm{S^n}{\pi}\xhookrightarrow{} \avrm{S^n}{r}$.
Thus for $0<r\leq\pi$, there exists an embedded copy of $\RP^n$ inside $\avrm{S^n}{r}$ given by the map $[x]\mapsto \frac{1}{2}\delta_x + \frac{1}{2} \delta_{-x}$.
We are going to refer this copy of $\RP^n$ inside of $\avrm{S^n}{r}$ as the \emph{central core}.

\begin{lemma}
\label{lem:piCase}
$\avrm{S^n}{\pi}\simeq \RP^n$ for all $n\geq 1$.
\end{lemma}

\begin{proof}
We will construct a deformation retration of $\avrm{S^n}{\pi}$ onto the central core $\RP^n$.
Define $H\colon \avrm{S^n}{\pi}\times I\to \avrm{S^n}{\pi}$ as follows.
For $\mu = \lambda\delta_x + (1-\lambda)\delta_{-x}\in \avrm{S^n}{\pi}$ with $0\leq \lambda \leq 1$, define
\[H(\mu,t) = (1-t)\mu + t(\tfrac{1}{2}\delta_x + \tfrac{1}{2}\delta_{-x})= \left(\lambda(1-t)+\tfrac{t}{2}\right)\delta_x+\left((1-\lambda)(1-t)+\tfrac{t}{2}\right)\delta_{-x}.\]
Note that for $\nu = \frac{1}{2}\delta_x + \frac{1}{2}\delta_{-x}\in \phi(\RP^n)$, we have $H(\nu,t) = \nu$ for all $t\in I$.
Moreover, for $\mu =\lambda\delta_x + (1-\lambda)\delta_{-x}$, we have $H(\mu,0) = \mu$ and $H(\mu,1) = \frac{1}{2}\delta_x + \frac{1}{2}\delta_{-x}\in \phi(\RP^n)$.
Since $H$ is well-defined and since $\mu\mapsto \tfrac{1}{2}\delta_x + \tfrac{1}{2}\delta_{-x}$ is continuous, the map $H$ is continuous by Lemma~\ref{lem:3.9ofAAR}.
Therfore, $H$ is a deformation retraction of $\avrm{S^n}{\pi}$ onto $\phi(\RP^n) \cong \RP^n$, completing the proof.

\end{proof}

\begin{remark}
We will refer to the map $H(\,\cdot\,, 1)\colon \avrm{S^n}{\pi} \to \avrm{S^n}{\pi}$ given by $\lambda \delta_x+(1-\lambda)\delta_{-x} \mapsto \frac{1}{2}(\delta_x + \delta_{-x})$ %$, \forall 0\leq \lambda \leq 1, \forall x\in S^n$ 
as \emph{projection to the central core}.
\end{remark}

In Subsection~\ref{subsec:fiber_bundle}, we describe a fiber bundle structure on $\avrm{S^n}{\pi}$.
Using fibrations and Whitehead's theorem, we also give an alternate proof of Lemma~\ref{lem:piCase}.

We are now prepared to prove Theorem~\ref{thm:avrmSn-homotopy-type}.

\begin{proof}[Proof of Theorem~\ref{thm:avrmSn-homotopy-type}]

We must show that $\avrm{S^n}{r}\simeq \RP^n$ for all $\frac{2\pi}{3}<r\le\pi$.
We will accomplish this by showing that $\avrm{S^n}{r}$ is homotopy equivalent to $\avrm{S^n}{\pi}$, which is homotopy equivalent to $\RP^n$ by Lemma~\ref{lem:piCase}.

We have the canonical inclusion $\iota\colon \avrm{S^n}{\pi}\xhookrightarrow{} \avrm{S^n}{r}$.
We define a homotopy inverse for $\iota$, denoted $\rho\colon \avrm{S^n}{r} \to \avrm{S^n}{\pi}$, as follows.
First, for measures supported on antipodal points, we define $\rho(\lambda_0 \delta_{x_0} + (1-\lambda_0)\delta_{-x_0})=\lambda_0 \delta_{x_0} + (1-\lambda_0)\delta_{-x_0}$ for all $x_0\in S^n$ and $0\leq \lambda_0\leq 1$.
In general, let $\mu=\lambda_0 \delta_{x_0} + (1-\lambda_0)\delta_{x'_0}$, where $x_0,x'_0\in S^n$, and where $r\leq d(x_0,x'_0)\leq\pi$.
Since $S^n$ sits inside the (convex) Euclidean space $\R^{n+1}$, we can consider the associated point $\tilde{\mu} = \lambda_0 x_0 + (1-\lambda_0)x'_0\in \R^{n+1}$.
%Note that so long as $x\neq x'$ and $0<\lambda<1$, then we have $\tilde{\mu}\in \R^{n+1}\setminus S^n$.
Since $\frac{2\pi}{3}<r\leq d(x_0,x'_0)\leq\pi$, there is a unique\footnote{
When $x_0'=-x_0$, meaning $d(x_0,x_0')=\pi$, there is no \emph{unique} great circle passing through $x_0$ and $x^\prime_0$, in fact, there are infinitely many.
But the image $\rho(\mu)$ is not dependent on which great circle is chosen.}
great circle in $S^n$ passing through $x_0$ and $x^\prime_0$.
Let $N$ be the ``north pole'' point on that great circle, by which we mean that $N$ is the point on the great circle such that $d(N,x_0)=d(N,x_0^\prime)=\frac{d(x_0,x_0^\prime)}{2}$.
Also, let $x_1$ and $-x_1$ be the points on the great circle such that $d(N,x_1)=d(N,-x_1)=\frac{\pi}{2}$;
note that the line in $\R^{n+1}$ between $x_0$ and $x'_0$ is parallel to the line in $\R^{n+1}$ between $x_1$ and $-x_1$.
We join $N$ with $\tilde{\mu}$ via a straight line in $\R^{n+1}$.
If the line intersects a point on the diameter 
$\{tx_1+(1-t)(-x_1)\in\R^{n+1}~|~0\le t\le 1\}$, then let the intersection point be $\lambda_1 x_1+(1-\lambda_1)(-x_1)$, and define $\rho(\mu)=\lambda_1\delta_{x_1}+(1-\lambda_1)\delta_{-x_1}$.
Otherwise, let $v$ be the unique point where this line intersects $S^n\setminus N$, and then define $\rho(\mu)=\delta_v$.
Even though it is defined in piecewise manner, $\rho\colon \avrm{S^n}{r} \to \avrm{S^n}{\pi}$ is a continuous function.
The map $\rho$ is depicted in Figure~\ref{fig:flashlight_homotopy} (left).

\begin{figure}[htb]
\def\svgwidth{3in}
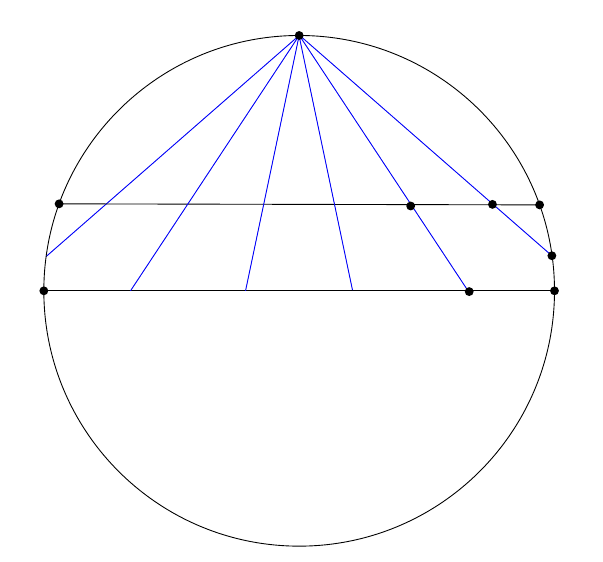
\hspace{3mm}
\def\svgwidth{3.2in}
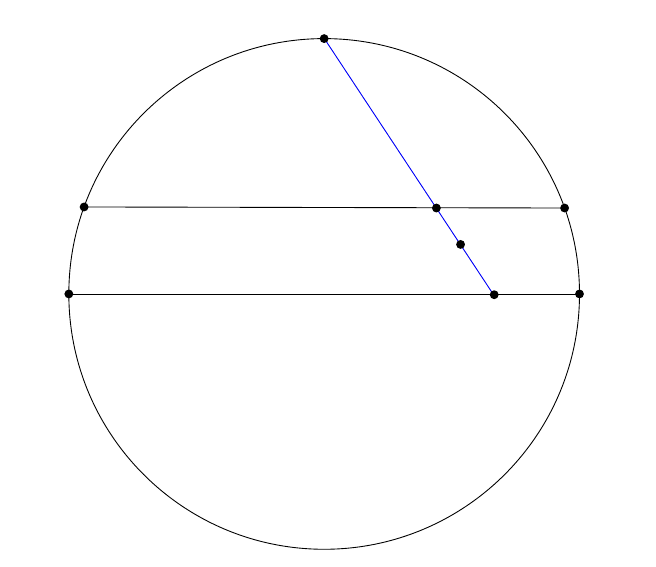
\caption{
(\emph{Left})
The map $\rho\colon \avrm{S^n}{r} \to \avrm{S^n}{\pi}$ in the proof of Theorem~\ref{thm:avrmSn-homotopy-type}, with $n=1$.
There are the two cases for how $\rho$ is defined.
In one case, $\mu = \lambda_0\delta_{x_0}+(1-\lambda_0)\delta_{x'_0}$ and $\rho(\mu)=\lambda_1\delta_{x_1}+(1-\lambda_1)\delta_{-x_1}$.
In the other case, $\rho(\mu')=\delta_{v}$.
(\emph{Right}) A depiction of the flashlight homotopy $H$ between $id_{\avrm{S^n}{r}}$ and $\rho$.
Here $H(\mu,t)=\lambda_t \delta_{x_{t}}+(1-\lambda_t)\delta_{x'_{t}}$.
}
\label{fig:flashlight_homotopy}
\end{figure}

We claim the following:
\begin{itemize}
\item $\rho \circ \iota=id_{\avrm{S^n}{\pi}}$ (trivial from the definition)
\item $\iota \circ \rho \simeq id_{\avrm{S^n}{r}}$.
\end{itemize}

To prove the second bullet point, consider the homotopy $H\colon \avrm{S^n}{r}\times I \to \avrm{S^n}{r}$, defined as follows.
Let $\widetilde{\rho}(\mu)\in \R^{n+1}$ be the point in Euclidean space corresponding to the measure $\rho(\mu)\in \avrm{S^n}{\pi}$, meaning either $\widetilde{\rho}(\mu)=\lambda_1 x_1+(1-\lambda_1)(-x_1)$ or $\widetilde{\rho}(\mu)=v$, according to the two cases in the definition of $\rho$.
For $t\in I$, let $\widetilde{H}(\mu,t)\in \R^{n+1}$ be defined by $\widetilde{H}(\mu,t)=t\tilde{\mu}+(1-t)\widetilde{\rho}(\mu)$.
Note that $\widetilde{H}(\mu,t)$ can be uniquely expressed as a convex combination 
$\widetilde{H}(\mu,t)=\lambda_t x_{t,\mu}+(1-\lambda_t)x'_{t,\mu}$ where $x_{t,\mu}\in S^n$ is a point on the geodesic arc from $x_0$ to $x_1$, where $x'_{t,\mu}$ is a point on the geodesic arc from $x'_0$ to $-x_1$, and where the line segment from $x_{t,\mu}$ to $x'_{t,\mu}$ is parallel to the line segment from $x_0$ to $x'_0$.
We emphasize that $x_{t,\mu}$ and $x'_{t,\mu}$ are dependent on both $t$ and $\mu$, which justifies the subscripts.
We note $x_0=x_{0,\mu}$, $x'_0=x'_{0,\mu}$, $x_1=x_{1,\mu}$, and $-x_1=x'_{1,\mu}$.
And, when it is clear from context, we may write $x_{t,\mu}$ simply as $x_t$, and $x'_{t,\mu}$ as $x'_t$, as we do for example in Figure~\ref{fig:flashlight_homotopy}.
We define $H(\mu,t)=\lambda_t \delta_{x_{t,\mu}}+(1-\lambda_t)\delta_{x'_{t,\mu}}$.
We refer to this as the \emph{flashlight homotopy}; see Figure~\ref{fig:flashlight_homotopy} (right).
We note that the homotopy $H$ satisfies $H(\,\cdot\,,0)=id_{\avrm{S^n}{r}}$ and $H(\,\cdot\,,1)=\iota \circ \rho$, as desired.

We are going to show that $H$ is a continuous function.
Since $\avrm{S^n}{r}\times I$ is a metric space it suffices to show that for every convergent sequence $(\mu_n,t_n)\to (\mu,t)$, we have a convergent sequence $H(\mu_n,t_n)\to H(\mu,t)$.
We can use the metric defined via $d\big((\mu_1,t_1),(\mu_2,t_2)\big)=\max (W_1(\mu_1,\mu_2),|t_1-t_2|).$
Also, we note that $W_1$ metrizes weak convergence~\cite[Theorem~6.9]{VillaniOldNew}.
This means that $\mu_n\xrightarrow{w}\mu$ ($\mu_n$ converges weakly to $\mu)$ , i.e., $\int fd\mu_n\to \int fd\mu$ for all continuous functions $f\colon S^n\to \R$, if and only if $W_1(\mu_n,\mu)\to 0$.

Suppose $(\mu_n,t_n)\to (\mu, t)$.
We want to show that $H(\mu_n,t_n)\to H(\mu,t)$.
Since $(\mu_n,t_n)\to (\mu, t)$, we have $|t_n-t|\to 0$ and $W_1(\mu_n,\mu)\to 0$, and thus $\int fd\mu_n\to \int fd\mu$ for all continuous functions $f\colon S^n\to \R$.
Now 
\[ W_1(H(\mu_n,t_n),H(\mu,t))=W_1(\underbrace{\lambda_{t_n}\delta_{x_{t_n,\mu_n}}+(1-\lambda_{t_n})\delta_{x^\prime_{t_n,\mu_n}}}_{\nu_{t_n}},\underbrace{\lambda_t\delta_{x_{t,\mu}}+(1-\lambda_t)\delta_{x^\prime_{t,\mu}}}_{\nu_{t}}).\]
Let $f\colon S^1\to \R$ be a continuous function.
Note
\begin{align*}
\int fd\nu_{t_n}&=\lambda_{t_n}f(x_{t_n,\mu_n})+(1-\lambda_{t_n})f(x^\prime_{t_n,\mu_n})\\
&\to \lambda_tf(x_{t,\mu})+(1-\lambda_{t})f(x^\prime_{t,\mu})=\int fd\nu
\end{align*}
\noindent since $\lambda_{t_n}\to \lambda_t$, since $x_{t_n,\mu_n}\to x_{t,\mu}$, and since $f$ is continuous.
Hence $H(\mu_n,t_n)\to H(\mu,t)$, and so $H$ is continuous.
\end{proof}

We can do even better.
At scale $r=\frac{2\pi}{3}$, if we ignore all the measures that are supported on 3 points, the resulting space $\avrm{S^n}{\frac{2\pi}{3}}\setminus W$, where $W$ is the set of all measures supported on 3 points, is homotopy equivalent to $\RP^n$.
We will use the following lemma in Theorem~\ref{thm:hom_type_avrm(S^n;2pi/3)}.

\begin{lemma}
\label{lem:2pi/3 equal case}   
Let $W=\{\lambda_1\delta_{x_1}+ \lambda_2\delta_{x_3}+\lambda_1\delta_{x_3}\in \avrm{S^n}{\frac{2\pi}{3}}\mid \lambda_i>0\ \forall \ i\}$.
Then $\avrm{S^n}{\frac{2\pi}{3}}\setminus W\simeq \RP^n$.
\end{lemma}

Note that $\avrm{S^n}{\frac{2\pi}{3}}$ has measures supported on at most $3$ points.
If a measure in $\avrm{S^n}{\frac{2\pi}{3}}$ is supported on $3$ points, the support is the set of vertices of an equilateral triangle.
The set of these measures is exactly the set $W$ in Lemma~\ref{lem:2pi/3 equal case}.\

\begin{proof}[Proof of Lemma~\ref{lem:2pi/3 equal case}]
We have natural inclusion $\iota \colon \avrm{S^n}{\pi} \xhookrightarrow{} \avrm{S^n}{\frac{2\pi}{3}}\setminus W$.
Again, the homotopy inverse of $\iota$ is $\rho \colon \avrm{S^n}{\frac{2\pi}{3}}\setminus W \to \avrm{S^n}{\pi}$, as defined in the proof of Theorem~\ref{thm:avrmSn-homotopy-type}.
The map $\rho$ can be defined the same way, since all measures in $\avrm{S^n}{\frac{2\pi}{3}}\setminus W$ have support on at most two points.
\end{proof}

\subsection{A fiber bundle perspective on anti-VR thickenings}
\label{subsec:fiber_bundle}

In this subsection we give a fiber bundle structure of $\avrm{S^n}{\pi}$.
Note that $\avrm{S^n}{\pi}$ is a fiber bundle over $\RP^n$ with fibers homeomorphic to $I$.
Define $p \colon \avrm{S^n}{\pi}\to \RP^n$ by $t\delta_x+(1-t)\delta_{-x}\mapsto [x]$.
Note $p^{-1}([x])=\{t\delta_x+(1-t)\delta_{-x}\mid t\in I\}\cong I$.
To see that this is a fiber bundle, consider $[x]\in \RP^n$.
Let $x\in S^n$, and let $V$ be a neighborhood of $x$ in $S^n$ that does not contain any antipodal pair.
Let $U$ be the image of $V$ under the antipodal identification.
Then $p^{-1}(U)$ consists of all of the points $t\delta_x+(1-t)\delta_{-x}$ with $[x]\in U.$
Define $h:p^{-1}(U)\to U\times I$ by $t\delta_x+(1-t)\delta_{-x}\mapsto ([x],t)$.
Now $h$ is a continuous map with continuous inverse $([x],t)\mapsto t\delta_x+(1-t)\delta_{-x}$.

We note a similarity with the tautological line bundle $\gamma_n$ on projective space $\RP^n$.
Recall that the total space of $\gamma_n$ consists of the points $\{([x],tx)\in \RP^n\times \R^{n+1}\mid t\in \R\}$.
That means the fibers are homeomorphic to $\R$, but in the case of $p\colon \avrm{S^n}{\pi}\to \RP^n$, the fibers are instead homeomorphic to the unit interval.

Also, if we compose the map $S^n \to \avrm{S^n}{\pi}$ defined by $x\mapsto \delta_x$ with the map $p\colon \avrm{S^n}{\pi} \to \RP^n$, then we obtain the 2-fold universal covering map $S^n\to \RP^n$ defined via $x \mapsto [x]$.

Note the following facts:
\begin{itemize}
\item Every fiber bundle is Serre fibration.

\item  Let $p \colon E\to B$ be a (Hurewicz) fibration, with base space $B$ path connected.
If $B$ and $F=p^{-1}(b)$ have the homotopy type of a CW complex for every $b\in B$, then $E$ has the homotopy type of a CW complex~\cite[Theorem~5.4.2]{Fritsch_Piccinini_1990}.

\item For every Serre fibration $p \colon E \to B$, there exists a long exact sequence in homotopy groups.
For $b\in B$, $x\in p^{-1}(b)=F$ we have 
\[\ldots \to \pi_n(F,x)\to \pi_n(E,x)\to \pi_n(B,b)\to \pi_{n-1}(F,x)\to \ldots,\]
where $\pi_n(F,x)\to \pi_n(E,x)$ is induced by the inclusion map and $\pi_n(E,x)\to \pi_n(B,b)$ is induced by $p$.
\end{itemize}

In the case of $p \colon \avrm{S^n}{\pi}\to \RP^n$, the fibers are contractible.
Using the second fact, we conclude that $\avrm{S^n}{\pi}$ has the homotopy type of a CW complex.
Using the first and third fact, we get that $\pi_n(\avrm{S^n}{\pi})\cong \pi_n(\RP^n)$ for all $n$.
From Whitehead's theorem, we conclude that $\avrm{S^n}{\pi}\simeq \RP^n$, giving another proof of Lemma~\ref{lem:piCase}.

See Section~\ref{sec:no-graph-homomorphisms} for yet another proof of Lemma~\ref{lem:piCase}, using properties of mapping cylinders.

\section{The covering dimension of the anti-VR thickening of a compact manifold}
\label{sec:cov-dim}

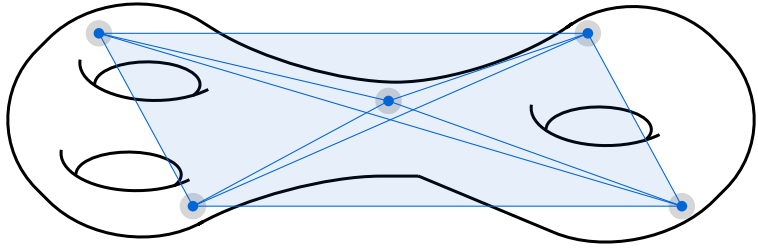
\begin{figure}[htb]
\begin{center}
\begin{tikzpicture}

\draw[line width=.4mm, rounded corners = 37](6,-.75)--(4.2,-.75)--(2,-2)--(0,0)--(2,2)--(4.4,.5)--(7,.5)--(9,2)--(11,0)--(9,-2)--(6,-.75);
\draw[line width=.4mm,] (1.25,-.4) arc (175:315:1cm and 0.5cm);
\draw[line width=.4mm,] (2.75,-0.88) arc (-30:180:0.7cm and 0.3cm);
%\draw (5.8,0) arc (0:360:0.5cm and 1cm);
\draw[line width=.4mm,] (1.5,.8) arc (175:315:1cm and 0.5cm);
\draw[line width=.4mm,] (3,.32) arc (-30:180:0.7cm and 0.3cm);
\draw[line width=.4mm,] (7.5,0.2) arc (175:315:1cm and 0.5cm);
\draw[line width=.4mm,] (9,-0.28) arc (-30:180:0.7cm and 0.3cm);
%\node[thick] (a) at (-13:5.8) {$\partial M$};
%\node[thick] (a) at (26:2.5) {$\tilde{M}$};

% Define points
\coordinate (A) at (1.75,1.15); % Left-top
\coordinate (B) at (8.25,1.15);   % Right-top
\coordinate (C) at (3.0,-1.15);% Left-bottom
\coordinate (D) at (9.5,-1.15);  % Right-bottom
\coordinate (E) at (5.6,.25);  % Center
		
% Add gray highlights around nodes
\foreach \point in {A,B,C,D,E} {
    \fill[mygray] (\point) circle (5pt);
}
% Draw nodes
\foreach \point in {A,B,C,D,E} {
    \fill[myblue] (\point) circle (2pt);
}

% Draw shaded area
\fill[myblue, opacity=0.1] (A) -- (B) -- (D) -- (C) -- cycle;
		
% Draw edges
\draw[myblue] (A) -- (B);
\draw[myblue] (B) -- (D);
\draw[myblue] (D) -- (C);
\draw[myblue] (C) -- (A);
\draw[myblue] (A) -- (E);
\draw[myblue] (B) -- (E);
\draw[myblue] (C) -- (E);
\draw[myblue] (D) -- (E);
\draw[myblue] (A) -- (D);
\draw[myblue] (C) -- (B);

\end{tikzpicture}
\end{center}
\caption{For $M$ a compact $n$-dimensional manifold with $p=\pack_M(r)$, the covering dimension of $\avrm{M}{r}$ is at most $(n+1)p-1=(p-1)+np$.
Above is a manifold $M$ of dimension $n=2$ with $r>0$ such that $p=\pack_M(r)=5$.
In blue is a single $(p-1)$-dimensional simplex in $\avrm{M}{r}$.
The gray open neighborhoods show that each of the $p$ vertices of this simplex have $n$ degrees of freedom.
Hence we expect $\avrm{M}{r}$ to have dimension $(p-1)+np=4+2\cdot 5=14$.}
\label{fig:cov-dim}
\end{figure}

Let $M$ be a compact $n$-dimensional Riemannian manifold.
Recall the \emph{$r$-packing number} $p\coloneqq\pack_M(r)$ is the maximal number of points that can be placed in $M$ such that any two points are at distance at least $r$.
So any probability measure in $\avrm{M}{r}$ has support of size at most $p$.
Is there some sense in which $\avrm{M}{r}$ has dimension at most $(n+1)p-1=(p-1)+np$?
This back-of-the-envelope guess at a numeric bound arises since a simplex on $p$ points is $(p-1)$-dimensional, and if the vertices of a simplex live in an $n$-dimensional manifold then each of the $p$ vertices has $n$ degrees of freedom; see Figure~\ref{fig:cov-dim}.

In Theorem~\ref{thm:cov-dim} we show the answer is yes: the \emph{covering dimension} of $\avrm{M}{r}$ is at most $(n+1)p-1$.
As a consequence, in Corollary~\ref{cor:cov-dim} we conclude that the cohomology $H^k(\avrm{M}{r})$ vanishes for all $k > (n+1)p-1$.

The \emph{covering dimension} of a topological space $X$ is at most $m$ if for every open cover $\cU$ of $X$, there is an open refinement $\cU'\subseteq \cU$ of multiplicity at most $m+1$; see~\cite[Section~50]{munkres} and~\cite[Page~152]{spanier1989algebraic}.
The refinement $\cU'$ is said to have multiplicity at most $m+1$ if no point in $X$ is contained in more than $m+1$ sets in $\cU'$.
% https://en.wikipedia.org/wiki/Lebesgue_covering_dimension

We use an equivalent definition of the covering number for compact metric spaces.
For $X$ a compact metric space, the \emph{covering dimension} of $X$ is at most $m$ if for any $\varepsilon > 0$, there is a cover of multiplicity at most $m+1$ by open sets of diameter at most $\varepsilon$.
% For the equivalence, the proof of => is easy, and the proof of <= should use Lebesgue's number lemma https://en.wikipedia.org/wiki/Lebesgue%27s_number_lemma.
Indeed, this follows from~\cite[Theorem~2.2(2)]{DRANISHNIKOV2018429}, which instead presents a cover of multiplicity at most $m+1$ as an $(m+1)$-fold union $\cU^0\cup\ldots\cup\cU^m$ where each $\cU^i$ is a family of disjoint sets.

The following result is related to~\cite[(3.1)]{fedorchuk1991probability} and~\cite{basmanov1983covariant}.

\begin{theorem}
\label{thm:cov-dim}
Let $M$ be a compact Riemannian manifold of dimension $n$, with $p=\pack_M(r)$.
Then the covering dimension of $\avrm{M}{r}$ is at most $(n+1)p-1$.
\end{theorem}

\begin{proof}
Let $\varepsilon>0$.
Since $M$ is compact, by Lemma~\ref{lem:compact} we know that $\avrm{M}{r}$ is compact.
Therefore, it suffices to show there exists a cover of $\avrm{M}{r}$ of multiplicity at most $(n+1)p$ by open sets of diameter at most $\varepsilon$.

We know that the compact $n$-dimensional manifold $M$ has covering dimension $n$~\cite[Corollary~50.7]{munkres}.
Therefore, there exists a cover $\cU$ of $M$ of multiplicity at most $n+1$ by open sets of diameter less than $\frac{\varepsilon}{2}$.
Since $M$ is compact, we can take this open cover $\cU$ to be finite.

It suffices to restrict attention to $\varepsilon$ satisfying $0<\varepsilon<r$, which we do.
It therefore follows that given any $\mu\in\avrm{M}{r}$ and $U\in \cU$, the intersection $\supp(\mu)\cap U$ consists of at most one point.

Let $K$ be the simplicial complex with vertex set $\cU$ whose simplices are the subsets $\{U_0,\ldots,U_k\}\subseteq \cU$ such that there is some $\mu \in \avrm{M}{r}$ with $\mu(U_i)>0$ for all $0\le i\le k$.
Since each $\mu \in \avrm{M}{r}$ has support of size at most $p$, and since each point in $M$ is contained in at most $n+1$ sets in $\cU$, each simplex in $K$ has at most $(n+1)p$ vertices.
Hence the simplicial complex $K$ is of dimension at most $(n+1)p-1$.

Any point $w\in |K|$ in the geometric realization of $K$ is a collection of nonnegative weights $\{w_U\}_{U\in \cU}$, at most $(n+1)p$ of them positive, whose sum is $1$.
We might refer to $w\in |K|$ as a \emph{weight sequence}.
Equip the geometric realization $|K|$ with the metric of barycentric coordinates (\cite[Section~7A.5]{bridson2011metric} and~\cite{dowker1952topology}).
In general, the metric of barycentric coordinates only preserves the weak homotopy type of a simplicial complex, but here it preserves the space up to homeomorphism since $K$ is finite.
Let $\varepsilon'>0$.
Since the geometric realization of a simplicial complex of dimension at most $(n+1)p-1$ has covering dimension at most $(n+1)p-1$~\cite[Page~152]{spanier1989algebraic}, 
there exists an open cover $\cW$ of $|K|$ by sets of diameter at most $\varepsilon'$ such that each point of $|K|$ belongs to at most $(n+1)p$ sets in the cover.
We choose $\varepsilon'<\frac{\varepsilon}{2\diam(M)}$, where $\diam(M)<\infty$ since $M$ is compact.
Unravelling the definition of the metric of barycentric coordinates, this means that within each set $W\in \cW$ in the cover, the sum of weight differences at all of the vertices $U\in \cU$ is less than $\varepsilon'$.
In symbols, for any $w,w'\in W\in \cW$, we have $\sum_{U\in \cU}|w_U-w'_U|<\varepsilon'$.

Let $\{\rho_U\}_{U\in \cU}$ be a partition of unity subordinate to $\cU$, which means that each $\rho_U \colon M \to \R_{\ge 0}$ satisfies $\supp(\rho_U)\subseteq U$, and each $x\in M$ satisfies $\sum_{U\in \cU}\rho_U(x)=1$.
For a measure $\mu \in \avrm{M}{r}$ we have that $U\cap\supp(\mu)$ is either a single point or the emptyset, since $\varepsilon<r$.
Therefore, any measure $\mu \in \avrm{M}{r}$ can be written as 
\[\mu 
= \sum_{\{U\in \cU:\mu(U)>0\}} \mu(U)\rho_U(x_U)\delta_{x_U} 
= \sum_{U\in \cU} \mu(U)\rho_U(x_U)\delta_{x_U},\]
where $x_U$ is the unique point in $U\cap\supp(\mu)$ for all $U\in\cU$ with $\mu(U)>0$ (and where the location of $x_U$ is irrelevant if $\mu(U)=0$).

We define a continuous map $h\colon \avrm{M}{r}\to |K|$ using this partition of unity subordinate to the cover $\cU$.
Indeed, let $\mu=\sum_{U\in \cU}\mu(U)\rho_U(x_U)\delta_{x_U}\in \avrm{M}{r}$ map to the weight sequence $\{h(\mu)_U\}_{U\in \cU}\in |K|$ defined by
$
h(\mu)_U = \mu(U)\rho_U(x_U)
$,
where again the location of $x_U$ is irrelevant if $\mu(U)=0$.
Note that $\{h(\mu)_U\}_{U\in \cU}$ is a weight sequence since $\mu$ is a probability measure.

Now, define the open cover $\cV$ of $\avrm{M}{r}$ by $\cV=\{h^{-1}(W)\}_{W\in \cW}$.
Since $\cW$ covers $|K|$, we have that $\cV$ covers $\avrm{M}{r}$.
Since the sets $W\in \cW$ are open in $|K|$ and since $h$ is continuous, the sets $h^{-1}(W)$ are open in $\avrm{M}{r}$.
So $\cV$ is an \emph{open} cover of $\avrm{M}{r}$.
We furthermore claim that the multiplicity of $\cV$ is at most the multiplicity of $\cW$, which is at most $(n+1)p$.
Indeed, note that if the sets $h^{-1}(W_1)$, \ldots, $h^{-1}(W_k) \in \cV$ have nonempty intersection at a point $\mu \in \avrm{M}{r}$, then the sets $W_1,\ldots,W_k$ have nonempty intersection at the point $h(\mu)\in |K|$.

We now show how to bound the diameter of a set $V_W \in \cV$.
Let $\mu=\sum_{U\in \cU}\mu(U)\rho_U(x_U)\delta_{x_U}$ be an arbitrary measure in $\avrm{M}{r}$, where $x_U$ is the unique point in $U\cap\supp(\mu)$ for all $U\in\cU$ with $\mu(U)>0$.
Similarly, let $\mu'=\sum_{U\in \cU}\mu'(U)\rho_U(x'_U)\delta_{x'_U}\in \avrm{M}{r}$, where $x'_U\in U$ is defined similarly.
Suppose $\mu,\mu'\in V_W$ for some $V_W\in \cV$.
So $h(\mu),h(\mu')\in W$ for some $W\in \cW$.
Since $W$ has diameter at most $\varepsilon'$, this means that
\[
\varepsilon' > \sum_{U\in \cU}|h(\mu)_U-h(\mu')_U| 
%&= \sum_{U\in \cU}|\mu(U)\rho_U(U\cap \supp(\mu))-\mu'(U)\rho_U(U\cap \supp(\mu'))| \\
%&\ge \sum_{\{U:\mu(U)\ge\mu'(U)\}}|\mu(U)\rho_U(U\cap \supp(\mu))-\mu'(U)\rho_U(U\cap \supp(\mu'))| + \sum_{\{U:\mu'(U)>\mu(U)\}}\ldots\\
= \sum_{U\in \cU}|\mu(U)\rho_U(x_U)-\mu'(U)\rho_U(x'_U)|.
\]
We construct a transport plan in $\cP(M \times M)$ with marginals $\mu$ and $\mu'$ as follows.
First, we define a \emph{partial transport plan} (which may have mass less than one) by $\sigma=\sum_{U\in \cU}\sigma_U\delta_{(x_U,x'_U)}$, where for $U\in \cU$ we set
\[
\sigma_U = \min\{\mu(U)\rho_U(x_U),\mu'(U)\rho_U(x'_U)\}.
\]

We note that the cost of the partial transport plan $\sigma$ satisfies
$\cost(\sigma) \le \frac{\varepsilon}{2}$,
since $\sigma$ only ever matches mass at points $x_U$ and $x'_U$ that lie in a common open set in $\cU$, and each such set has diameter less than $\frac{\varepsilon}{2}$.
The mass of the partial transport plan $\sigma$ satisfies
\begin{align*}
\mass(\sigma) 
&= \sum_{U\in \cU} \sigma_U
= \sum_{U\in \cU} \min\{\mu(U)\rho_U(x_U),\mu'(U)\rho_U(x'_U)\} \\
&\ge \sum_{U\in \cU}\left( \mu(U)\rho_U(x_U) - |\mu(U)\rho_U(x_U)-\mu'(U)\rho_U(x'_U)| \right)\\
&= \sum_{U\in \cU} \mu(U)\rho_U(x_U) - \sum_{U\in \cU} |\mu(U)\rho_U(x_U)-\mu'(U)\rho_U(x'_U)| \\
&> 1 - \varepsilon'
\end{align*}
Therefore by Lemma~\ref{lem:ExtnOfOptimalTransport}, $\sigma$ can be completed to a (full) transport plan $\tilde{\sigma} \in \cP(M \times M)$ with marginals $\mu$ and $\mu'$ whose cost is at most
\begin{align*}
\cost(\tilde{\sigma}) &\le \cost(\sigma)+(1-\mass(\sigma))\diam(M) \\
&\le \tfrac{\varepsilon}{2} + \varepsilon' \diam(M) \\
&\le \tfrac{\varepsilon}{2} + \frac{\varepsilon}{2\diam(M)} \diam(M) 
= \varepsilon.
\end{align*}
Hence the diameter of a set $V_W \in \cV$ is at most $\varepsilon$.
We have shown the covering dimension of $\avrm{M}{r}$ is at most $(n+1)p-1$.
%%%%%%%%%%
%(\footnote{\note{In the above proof, what did we end up using about $h$?
%The fact that it was continuous, and what else?
%Bounded distortion, Michael suggests?}})
%%%%%%%%%%
\end{proof}

\begin{corollary}
\label{cor:cov-dim}
Let $M$ be a compact Riemannian manifold of dimension $n$, with $p=\pack_M(r)$.
Then the %homology $H_k(\avrm{M}{r})$ and cohomology $H^k(\avrm{M}{r})$ 
\v{C}ech cohomology $\check{H}^k(\avrm{M}{r})$ 
vanishes in all dimensions $k\geq (n+1)p$.
\end{corollary}

\begin{proof}
It follows from the definition of \v{C}ech cohomology that if a topological space has covering dimension at most $d$, then that space has vanishing \v{C}ech cohomology in all dimensions $k>d$.
% This is stated for example at \url{https://mathoverflow.net/questions/43400/topological-dimension-versus-cohomological-dimension}; for example BCnrad writes ``The covering dimension clearly is an upper bound on non-vanishing for Cech cohomology''.
\end{proof}

\begin{example}
By Theorem~\ref{thm:cov-dim} the covering dimension of $\avrm{S^n}{r}$ is at most $(n+1)\pack_{S^n}(r)$, and by Corollary~\ref{cor:cov-dim} the cohomology of $\avrm{S^n}{r}$ vanishes in all dimensions $k\ge (n+1)\pack_{S^n}(r)$.
\end{example}

\begin{example}
Let us apply Theorem~\ref{thm:cov-dim} to $\avrm{S^1}{r}$ for $\frac{2\pi}{3}<r<\pi$.
We have that $S^1$ is a compact Riemannian manifold of dimension $n=1$, with $p=\pack_{S^1}(r)=2$ for $r$ in this range.
Therefore, Theorem~\ref{thm:cov-dim} states that the covering dimension of $\avrm{S^1}{r}$ is at most $(n+1)p-1=3$.
We note this bound is tight, since $\avrm{S^1}{r}$ for $\frac{2\pi}{3}<r<\pi$ is a $3$-dimensional space (which can be thought of as a 2-parameter family of edges) for $r$ in this range.
\end{example}

\begin{example}
However, for $r=\pi$, Theorem~\ref{thm:cov-dim} only implies that $\avrm{S^1}{\pi}$ has covering dimension at most $3$, even though (as a M\"obius band) it has covering dimension $2$.
This is because the points in $S^1$ which support edges in $\avrm{S^1}{\pi}$ are necessarily antipodal points, and hence can't be varied independently in small neighborhoods without decreasing the distance below $\pi$.
\end{example}

In the specific case of $\avrm{S^n}{\frac{2\pi}{3}}$, from Theorem~\ref{thm:cov-dim}, we conclude that the covering dimension will be at most $(n+1)\cdot 3 -1= 3n+2$.
Thus cohomology will vanish in dimension $3n+3$ and higher.
In the following section, we strengthen these bounds, showing the (co)homology of $\avrm{S^n}{\frac{2\pi}{3}}$ vanishes in dimensions $2n+2$ and higher.

\section{Anti-VR thickenings of spheres at scale $r = \frac{2 \pi}{3}$}
\label{sec:avr-small}

Recall that $\avrm{S^n}{r}$ is homeomorphic to $S^n$ for all $r > \pi$.
Also, from Theorem~\ref{thm:avrmSn-homotopy-type}, we know that $\avrm{S^n}{r}$ is homotopy equivalent to $\RP^n$ for all $\frac{2\pi}{3}<r\leq \pi$.
In this section we give some initial information about the topology of $\avr{S^n}{r}$ at the first scale, $r=\frac{2\pi}{3}$, when the homotopy type is no longer $\RP^n$.
Specifically, Theorem~\ref{thm:hom_type_avrm(S^n;2pi/3)} shows that $\avrm{S^n}{\frac{2\pi}{3}}$ has vanishing homology and cohomology in dimensions $2n+2$ and higher.

Intuitively, this result is obtained by recognizing pieces of $\avrm{S^n}{\frac{2\pi}{3}}$ of known dimensions as follows.
Up to homotopy, we start with $\RP^n$.
To this space, we glue on a $(2n-1)$-parameter family of solid triangles (disks), which is $(2n+1)$-dimensional.
The $2n-1$ parameter family arises since the planes in $\R^{n+1}$ are parametrized by $\gr(2,n+1)$, a manifold of dimension $2\cdot(n+1-2)=2n-2$, and within each plane we have a circle's worth of triangles, increasing the dimension by $1$.
We formalize this intuition in the remainder of this section.

\begin{theorem}
\label{thm:hom_type_avrm(S^n;2pi/3)}
%The space $\avrm{S^n}{\frac{2\pi}{3}}$ has the homotopy type of $(2n+1)$-dimensional CW complex.
We have vanishing cohomology $H^k(\avrm{S^n}{\tfrac{2 \pi}{3}}) \cong 0$ and vanishing homology $H_k(\avrm{S^n}{\tfrac{2 \pi}{3}}) \cong 0$ for all $n \geq 1$ and all $k \geq 2n+2$.
\end{theorem}

\begin{figure}[h]
\def\svgwidth{4in}
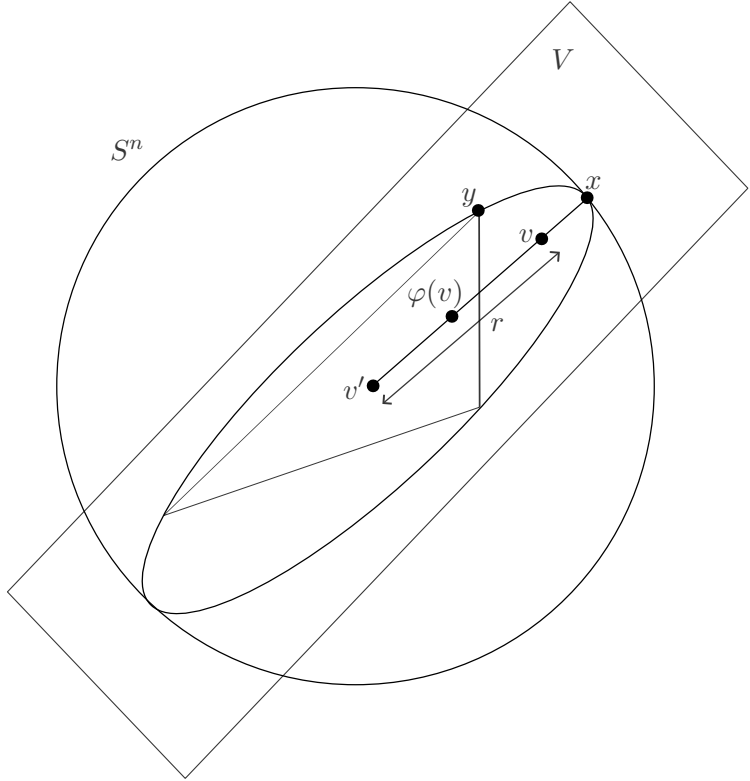
\caption{A \emph{pictorial} representation of an element $v$ of $Z \times I$ and its image under the map $\varphi \colon Z \times I \to \avrm{S^n}{\tfrac{2 \pi}{3}}$. 
Note that $v=(V,y,x,r)$ and $v'=(V,y,x,0)$ live in the space $Z\times I$ and $\varphi(v)$ lives in a different space $\avrm{S^n}{\frac{2\pi}{3}}$.
}
\label{fig:point in ZxI}
\end{figure}

\begin{proof}
We begin by introducing spaces that will allow us to parameterize the families of triangles described above:
\begin{align*}
Z &= \{ (V,y,x) \in \gr(2,n+1)\times S^n\times S^n~|~y,x\in V \} \\
C &= \{ (V,y) \in \gr(2,n+1)\times S^n~|~y\in V \}.
\end{align*}
Both are compact, as they are closed subsets of products of compact spaces (in fact, we will see below that both are closed manifolds).
We define a map $\varphi \colon Z \times I \to \avrm{S^n}{\tfrac{2 \pi}{3}}$ as follows.
For a given $(V,y,x,r) \in Z \times I$, there is a unique equilateral triangle centered at the origin that has $y$ as a vertex and is contained in $V$; let $y'$ and $y''$ be its other vertices.
Let $z$ be the point on the boundary of this triangle that lies on the line passing through the origin and $x$.
Then the point $rz$ can be written uniquely as a convex combination $ay+a'y'+a''y''$ where $a,a',a'' \geq 0$ and $a+a'+a''=1$; we define 
$\varphi(V,y,x,r) = a \delta_y + a' \delta_{y'} + a'' \delta_{y''}$, see Figure~\ref{fig:point in ZxI}.
In effect, this treats $x$ and $r$ as polar coordinates within the triangle specified by $V$ and $y$.

Since $Z \times I$ is compact and $\avrm{S^n}{\frac{2 \pi}{3}}$ is Hausdorff, $\varphi$ is a closed map and thus a quotient map onto its image.
The map $\varphi$ identifies points $(V,y,x,0)$ and $(V,y,x',0)$ for all $x$ and $x'$, and further identifies points $(V,y,x,t)$ and $(V,y',x,t)$ when the dot product is $\langle y , y' \rangle = -\frac{1}{2}$, i.e.\ when $y$ and $y'$ differ in angle by $\frac{2 \pi}{3}$.
Letting $W$ be defined as in Lemma~\ref{lem:2pi/3 equal case}, the image of $\varphi$ is the closure of $W$, and the points mapped by $\varphi$ into the closed set $\avrm{S^n}{\frac{2 \pi}{3}} \setminus W$ are exactly those of the form $(V,y,x,1)$.
Therefore, $\avrm{S^n}{\frac{2 \pi}{3}}$ is homeomorphic to the quotient
\[
\Big((Z\times I) \sqcup \avrm{S^n}{\tfrac{2 \pi}{3}} \setminus W\Big) / \sim,
\]
where $\sim$ is generated by $(V,y,x,t) \sim (V,y',x,t)$ if $\langle y , y' \rangle = -\frac{1}{2}$, $(V,y,x,0) \sim (V,y,x',0)$, and $(V,y,x,1) \sim \varphi(V,y,x,1)$.
This is equivalently a double mapping cylinder of a pair of maps 
\[ \widetilde{C} \xleftarrow[]{\psi} \widetilde{Z} \xrightarrow[]{\varphi_1} \avrm{S^n}{\tfrac{2 \pi}{3}} \setminus W.
\]
Here $\widetilde{C}$ and $\widetilde{Z}$ are quotients of $C$ and $Z$ formed by respectively identifying $(V,y)$ with $(V, y')$ and $(V,y,x) $ with $ (V, y', x)$ when $\langle y , y' \rangle = -\frac{1}{2}$, the map $\psi \colon \widetilde{Z} \to \widetilde{C}$ sends $(V,y,x)$ to $(V,y)$, and we define $\varphi_1(V,y,x) = \varphi(V,y,x,1)$.

We will check below that $\widetilde{C}$ and $\widetilde{Z}$ are closed, connected manifolds of dimensions $2n-1$ and $2n$ respectively; assuming this for now, we can get partial information about the homology of $\avrm{S^n}{\tfrac{2 \pi}{3}}$.
This double mapping cylinder above, homeomorphic to $\avrm{S^n}{\tfrac{2 \pi}{3}}$, can be written as the union of two mapping cylinders $M_\psi$ and $M_{\varphi_1}$ with intersection homeomorphic to $\widetilde{Z} \times I$, and we thus get Meyer--Vietoris sequences on homology and cohomology.
Since $\avrm{S^n}{\tfrac{2 \pi}{3}} \setminus W \simeq \RP^n$, we get the following long exact sequence on homology.
{\footnotesize
\[
% https://tikzcd.yichuanshen.de/#N4Igdg9gJgpgziAXAbVABwnAlgFyxMJZABgBpiBdUkANwEMAbAVxiRAB12oIcEBfUuky58hFAEZyVWoxZsAEgH0A1gApOAdyyw8DWMADCfAJQACThDTM4ppWs4BbOjgAWAIzfAASnwAKAPTBjEAEhbDwCIgAmKWp6ZlZEEDt1djoaACcHYABlQL5gThwAMwy6AGNgKPN2NCwCgGY+ExDBEAxw0SIG2JkEhUVgZQBacT5UrR0sPRhgAC0W0PbhCLFkABZe+LkkpSHR8c1tGF19IzMLKyYbPZGx1KdXD28-QOCljpFIlABWLdlEhwuDx+NIYFAAObwIigUoQBxIMggHAQJBjNpwhGISTI1GIKJLTFIGK4pBNDEZeFk6gopDrQmUrGbUmIH58Ch8IA
\begin{tikzcd}
\cdots \arrow[r] & H_k(\widetilde{C}) \oplus H_k(\RP^n) \arrow[r] & H_k(\avrm{S^n}{\tfrac{2 \pi}{3}}) \arrow[r] & H_{k-1}(\widetilde{Z}) \arrow[r] & H_{k-1}(\widetilde{C}) \oplus H_{k-1}(\RP^n) \arrow[r] & \cdots
\end{tikzcd}
\]
}
Since $\widetilde{C}$ has dimension $2n-1$, we have $H_k(\widetilde{C}) \oplus H_k(\RP^n) \cong 0$ for $k \geq 2n$, so $H_k(\avrm{S^n}{\tfrac{2 \pi}{3}}) \cong H_{k-1}(\widetilde{Z})$ for all $k \geq 2n+1$.
In particular, since $\widetilde{Z}$ is a manifold of dimension $2n$, this shows $H_k(\avrm{S^n}{\tfrac{2 \pi}{3}}) \cong 0$ for all $k \geq 2n+2$.
The analogous argument using the Meyer--Vietoris sequence for cohomology shows $H^k(\avrm{S^n}{\tfrac{2 \pi}{3}}) \cong 0$ for all $k \geq 2n+2$.

We now return to check that $\widetilde{C}$ and $\widetilde{Z}$ are closed, connected manifolds.
We begin with the spaces $C$ and $Z$, which we will in fact show are fiber bundles (see Section~\ref{sec:background} for background).
The space $C$ is defined similarly to the tautological vector bundle over $\gr(2, n+1)$, in which the fiber over $V \in \gr(2, n+1)$ is $V$ itself.
With this motivation, we check explicitly that the projection $\pi \colon C \to \gr(2, n+1)$ is a circle bundle.
We give $\mathbb{R}^{n+1}$  the usual inner product in order to work with orthogonal projections and orthogonal complements.
For a given $V \in \gr(2, n+1)$, let $P_V \colon \mathbb{R}^{n+1} \to V$ be the orthogonal projection, and let $U_V \subset \gr(2, n+1)$ be the open set of subspaces mapped bijectively onto $V$ by $P_V$.
Then for each $X \in U_V$, the map $X \cap S^n \to V \cap S^n$ given by $x \mapsto \frac{P_V(x)}{\|P_V(x)\|}$ is a homeomorphism, so we have a homeomorphism $\pi^{-1}(U_V) \cong U_V \times (V \cap S^n)$ that sends $(X,x)$ to $(X, \frac{P_V(x)}{\|P_V(x)\|})$.
Since $V \cap S^n \cong S^1$, this shows $C$ is a fiber bundle over $\gr(2, n+1)$ with fiber $S^1$.
An analogous argument shows that $Z$ is a fiber bundle over $\gr(2, n+1)$ with fiber $S^1 \times S^1$, and in particular, this implies $Z$ and $C$ are closed, connected manifolds of dimensions $2n$ and $2n-1$, respectively.
They further inherit smooth structures from the Grassmanian and the circle.

The quotient map $C \to \widetilde{C}$ identifies each point $(V,y)$ with $(V,y')$ and $(V,y'')$, where $y$, $y'$, and $y''$ are the vertices of an equilateral triangle lying in $V$.
We check that \emph{locally}, this map can be realized by quotienting by a group action that permutes the vertices of such triangles.
Beginning with the homeomorphism $\pi^{-1}(U_V) \cong U_V \times (V \cap S^n)$ above, fix an orientation of $(V \cap S^n) \cong S^1$ and assign each fiber $\pi^{-1}(X)$ the corresponding orientation under the homeomorphism.
For each $X \in U_V$, let $T_X$ be the linear transformation that rotates $X$ by an angle $\frac{2 \pi}{3}$ in the positively oriented direction and leaves the orthogonal complement of $X$ fixed.
For any $(X,x) \in \pi^{-1}(U_V)$, the vector $x^\perp$ that lies in $X$ and is at an angle of $\frac{\pi}{2}$ from $x$ in the positively oriented direction depends smoothly on $(X,x)$: it can be found by choosing a positively oriented orthonormal basis in $V$ containing $\frac{P_V(x)}{\|P_V(x)\|}$, applying the homeomorphism above to get a positively oriented basis in $X$ containing $x$, and applying the Gram--Schmidt process.
Since $T_X(x)$ can be written as a linear combination of $x$ and $x^\perp$, it is a smooth function of $(X,x)$.
We thus get a smooth action\footnote{
Note that this group action cannot be extended to all of $C$, since then it is not possible to assign consistent orientations to all fibers.
} of $\mathbb{Z}/3\mathbb{Z}$ on $\pi^{-1}(U_V)$ by setting $[n] \cdot (X,x) = (X, T_X^n x)$.
The quotient $\pi^{-1}(U_V) / (\mathbb{Z}/3\mathbb{Z})$ of $\pi^{-1}(U_V)$ by this group action identifies exactly those sets of points in $\pi^{-1}(U_V)$ that are identified by the quotient map $C \to \widetilde{C}$, and furthermore the quotient manifold theorem (\cite[Theorem~21.10]{lee2013smooth}) shows that $\pi^{-1}(U_V) / (\mathbb{Z}/3\mathbb{Z})$ is a manifold of dimension $2n-1$.
Therefore the class of any point $(V,y)$ in $\widetilde{C}$ has a neighborhood that is homeomorphic to $\mathbb{R}^{2n-1}$, so we conclude that $\widetilde{C}$ is a closed, connected manifold of dimension $2n-1$.
An analogous argument shows $\widetilde{Z}$ is a closed, connected manifold of dimension $2n$, as required in our work above, completing the proof of Theorem~\ref{thm:hom_type_avrm(S^n;2pi/3)}.
\end{proof}

We conjecture that the results of this section could be strengthened to show $\avrm{S^n}{\frac{2\pi}{3}}$ is in fact homotopy equivalent to a CW complex of dimension $2n+1$, since the fact that $\widetilde{C}$ and $\widetilde{Z}$ are manifolds imply they have CW structures (see~\cite[Corollary~A.12]{Hatcher} or~\cite[Theorem~3.5]{milnor1963morse}).

\section{The total anti-VR complex and thickening}
\label{sec:TAVR}

We now introduce the \emph{total anti-VR simpicial complex} and the \emph{total anti-VR metric thickening}, whose definitions and behavior are related to anti-VR complexes and thickenings.
The results in this section will be useful when proving Theorem~\ref{thm:no-graph-homomorphism} on the non-existence of graph homomorphisms between Borsuk graphs.
Indeed, the total anti-VR constructions allow more flexibility (than the anti-VR constructions) when working with the possibly discontinuous functions that arise in the context of graph homomorphisms.

\begin{definition}
\label{def:TAVR}
Let $X$ be a metric space and $r\ge 0$.
The \emph{total anti-VR complex} is the simplicial complex $\tavr{X}{r}$ with
\begin{itemize}
\item vertex set $X$
\item simplices of the form $\{x_1,...,x_k,y_1,\ldots,y_l\}$, where we require that $k,l\ge 1$, and that $d(x_i,y_j) \ge r$ for $1 \le i \le k$ and $1\le j\le l$ (and all faces of these simplices).
\end{itemize}
We let $\tavrm{X}{r}$ denote the associated \emph{total anti-VR metric thickening} from Definition~\ref{def:metric_thickening}.
\end{definition}

As we described in Section~\ref{ssec:anti-vrm}, in settings where the input metric space $X$ is not discrete, the topology of anti-VR metric thickenings does not seem to be so closely related to the topology of anti-VR simplicial complexes.
We expect, however, that the topology of the total anti-VR complex is more closely related to that of the anti-VR metric thickening.
Variants of the total anti-VR complexes which generalize from two sets $\{x_i\}_i$ and $\{y_j\}_j$ to more sets may be useful when considering anti-VR complexes and thickenings of spheres for scales below $\frac{2\pi}{3}$.

\begin{lemma}
\label{lem:atvr-inclusion}
Let $X$ be a metric space and $r \ge 0$.
Then we have inclusions $\avr{X}{r}\hookrightarrow\tavr{X}{r}$ and $\avrm{X}{r}\hookrightarrow\tavrm{X}{r}$.
\end{lemma}

\begin{proof}
It suffices to prove the simplicial complex case, which then implies the metric thickening case.
By definition, $X$ is the vertex set of $\tavr{X}{r}$, and it is also the vertex set of $\avr{X}{r}$.
Now, to see that every simplex of dimension at least one in $\avr{X}{r}$ is a simplex in $\tavr{X}{r}$, suppose $\{x_0,x_1,\ldots,x_k\}$ is a simplex in $\avr{X}{r}$ with $k\ge 1$.
This means that $d(x_i,x_j)\ge r$ for $0 \le i,j,\le k$ with $i\neq j$.
Hence if we let $y_1=x_0$, we see that $\{x_0,x_1\ldots,x_k\}=\{y_1,x_1\ldots,x_k\}$ is a simplex in $\tavr{X}{r}$, since $d(x_i,y_1)\ge r$ for all $1\le i\le k$.
\end{proof}

The following proposition will be used in the proof of the no-graph-homomorphism result in Theorem~\ref{thm:no-graph-homomorphism}.

\begin{proposition}
\label{prop:TAVR-to-AVR}
$\tavrm{S^n}{r}\simeq \RP^n$ for $\frac{2\pi}{3} < r  \le \pi$.
\end{proposition}

\begin{proof}
We will show that $\tavrm{S^n}{r}\simeq \avrm{S^n}{r}$ for all $r$ in this range, and then the result will follow since $\avrm{S^n}{r}\simeq \RP^n$ by Theorem~ \ref{thm:avrmSn-homotopy-type}.
Let $\mu = \sum_{i = 1}^k \mu_i\delta_{x_i} + \sum_{j = 1}^\ell \mu_j' \delta_{y_j}\in \tavrm{S^n}{r}$ with $\textup{supp}(\mu) = \{x_1,\dots,x_k, y_1, \dots, y_\ell\}$ satisfying $d(x_i,y_j)\geq r$ for all $1\leq i\leq k$ and $1\leq j\leq \ell$.

Let $\overline{\mu_x} = \sum_{i = 1}^k \tilde{\mu}_ix_i \in \R^{n+1}$ and $\overline{\mu_y} = \sum_{j = 1}^{\ell}\tilde{\mu}_j'y_j \in \R^{n+1}$, where $\tilde{\mu}_i = \mu_i/\sum_k \mu_k$ and $\tilde{\mu}_j' = \mu_j'/\sum_k \mu_k'$.
Consider the radial projections $\mu_x^* = \overline{\mu_x}/\norm{\overline{\mu_x}} \in S^n$ and $\mu_y^* = \overline{\mu_y}/\norm{\overline{\mu_y}} \in S^n$.
Define $\rho \colon \tavrm{S^n}{r} \to \avrm{S^n}{r}$ by $\mu\mapsto (\sum_{i = 1}^k\mu_i)\delta_{\mu_x^*} + (\sum_{j = 1}^\ell \mu_j')\delta_{\mu_y^*}$.
To see that $\rho$ is well-defined, we must show that $d(\mu_x^*,\mu_y^*) \ge r$.
Let $\theta_{ij} = d(x_i,y_j)$ for each $i$ and $j$.
Then $\theta_{ij} \ge r$ for $\frac{2\pi}{3} < r \le \pi$, giving that $\langle x_i, y_j\rangle = \cos(\theta_{ij}) \le \cos(r)$.
Hence, using the fact that $\sum_i\sum_j\tilde{\mu}_i\tilde{\mu}_j' = 1$, we have
\begin{align*}
\langle \overline{\mu_x}, \overline{\mu_y}\rangle
& = \textstyle \langle \sum_{i = 1}^k \tilde{\mu_i}x_i, \sum_{j = 1}^\ell \tilde{\mu}_j' y_j\rangle \textstyle
= \sum_i \sum_j \tilde{\mu}_i \tilde{\mu}_j' \langle x_i, y_j\rangle\\
& \textstyle \le \sum_i\sum_j \tilde{\mu}_i \tilde{\mu}_j'\cos(r)
= \cos(r)\left(\sum_i \tilde{\mu}_i\cdot \sum_j  \tilde{\mu}_j'\right)
= \cos(r)(1 \cdot 1)
= \cos(r).
\end{align*}
Since $\cos(r) < 0$ for $\frac{2\pi}{3} < r \le \pi$ and since $\norm{\overline{\mu_x}}\norm{\overline{\mu_y}}  \leq 1$, we thus have $\langle \overline{\mu_x}, \overline{\mu_y}\rangle \le \cos(r) \leq \norm{\overline{\mu_x}}\norm{\overline{\mu_y}} \cos(r)$.
It follows that $\langle \mu_x^*,\mu_y^*\rangle \le \cos(r)$ and hence $d(\mu_x^*,\mu_y^*) \ge r$, as desired.
Since $\mu_x^*$ and $\mu_y^*$ vary continuously with $\mu$, the map $\rho$ is continuous.

Let $\iota \colon \avrm{S^n}{r} \hookrightarrow \tavrm{S^n}{r}$ be the inclusion from Lemma~\ref{lem:atvr-inclusion}.
We claim that $\rho$ and $\iota$ are homotopy inverses.
Note that $\rho \circ \iota$ is the identity on $\avrm{S^n}{r}$.
For the other composition, we must show that $\iota \circ \rho$ is homotopic to the identity map on $\tavrm{S^n}{r}$.
To see this, consider the linear homotopy $H \colon \tavrm{S^n}{r} \times I \to \tavrm{S^n}{r}$ defined by
\begin{align*}
H(\mu,t) 
&= (1-t)\mu + t \cdot \iota(\rho(\mu)) \\
&= (1-t)\sum_{i = 1}^k \mu_i\delta_{x_i} + t \Bigg(\sum_{i = 1}^k\mu_i\Bigg)\delta_{\mu_x^*} + (1-t)\sum_{j = 1}^\ell \mu_j' \delta_{y_j} + t\Bigg(\sum_{j = 1}^\ell \mu_j'\Bigg)\delta_{\mu_y^*}.
\end{align*}
when $0<t<1$, the support of $H(\mu,t)$ is $\{x_1,...,x_k,\mu_x^*,y_1,\ldots,y_l,\mu_y^*\}$.
To see that $H$ is well-defined, we need to show $d(\mu_x^*,y_j) \ge r$ for $1\le j\le l$ and $d(x_i,\mu_y^*)\ge r$ for $1 \le i \le k$.

Using the fact that $\sum_{i = 1}^k \tilde{\mu_i}=1$, we get
\[
    \langle \overline{\mu_x}, y_j\rangle
    =\Big \langle \sum_{i = 1}^k \tilde{\mu_i}x_i,y_j\Big\rangle
    = \sum_{i = 1}^k \tilde{\mu_i} \langle x_i, y_j\rangle 
    \leq \sum_{i = 1}^k \tilde{\mu_i}\cos(r) 
    = \cos(r)\leq \|\overline{\mu_x}\|\cos(r).
\]
Thus, $\langle \mu_x^\ast, y_j\rangle\leq \cos(r)$, which proves that $d(\mu^\ast_x,y_j)\geq r$ for all $1\leq j\leq \ell$.
A similar argument shows $d(x_i,\mu_y^*)\ge r$.

The map $H$ is continuous by Lemma~\ref{lem:3.9ofAAR}.
Hence $H$ gives a homotopy from the identity map on $\tavrm{S^n}{r}$ to $\iota \circ \rho$.
We have shown $\tavrm{S^n}{r}\simeq \avrm{S^n}{r}$.
\end{proof}

We are interested in possible connections between total anti-Vietoris--Rips complexes, neighborhood complexes~\cite{Lovasz1978}, and box complexes~\cite{csorba2007homotopy,matousek2003using,matsushita2025dominance}.

\section{Chromatic numbers of Borsuk graphs}
\label{sec:choromatic-numbers-borsuk-grp}

This section contains a discussion of chromatic numbers in Section~\ref{ssec:chromatic}, Borsuk graphs in Section~\ref{ssec:borsuk}, and circular chromatic numbers and their relationship with Borsuk graphs in Section~\ref{ssec:circular-chromatic}.

\subsection{Chromatic numbers}
\label{ssec:chromatic}

Let $G=(V,E)$ be a simple graph, meaning that there are no loops and no multiple edges.

\begin{definition}
\label{def:graph-homomorphism}
A \emph{graph homomorphism} $f\colon G \to H$ from a graph $G$ to a graph $H$ is a function from the vertices of $G$ to the vertices of $H$ such that if $(u,v)$ is an edge in $G$, then $(f(u), f(v))$ is an edge in $H$.
\end{definition}

The \emph{chromatic number} of $G$, $\chi(G)$ is the minimum integer $k$ such that the vertex set $V$ can be colored with $k$ different colors, such that no two adjacent vertices have the same color.
Let $K_n$ be the complete graph on $n$ vertices.
Then in terms of graph homomorphisms,
\[\chi(G) = \inf\{n \mid \text{ there is a graph homomorphism } G \to K_n\}.\]
For example, $\chi(K_n)= n$.

\subsection{Borsuk graphs}
\label{ssec:borsuk}

Let $S^n$ be the $n$-sphere, equipped with the geodesic metric.
For any $0<\alpha<\pi$, let $G(S^n;\alpha)$ be the graph whose vertex set is the sphere $V=S^n$, and whose edge set is all pairs of vertices at distance exactly $\alpha$ apart.
% https://math.stackexchange.com/questions/4025794/definition-of-borsuk-graph
Similarly, let $\bor{S^n}{\alpha}$ be the graph whose vertex set is the sphere $V=S^n$, and whose edge set is all pairs of vertices at distance \emph{at least} $\alpha$ apart.
In different references, both the graph $G(S^n;\alpha)$ and the graph $\bor{S^n}{\alpha}$ have been given the name of the \emph{Borsuk graph}~\cite{Lovasz1978,lovasz1983self,raigorodskii2010chromatic,raigorodskii2012chromatic}.
The papers~\cite{kahle2020chromatic,prosanov2018chromatic} contain nice summaries as well.
The chromatic number of $G(S^n;\alpha)$ and of $\bor{S^n}{\alpha}$ now depend on the value of $\alpha$.

What are the chromatic numbers of $G(S^n;\alpha)$ and of $\bor{S^n}{\alpha}$ for $0<\alpha<\pi$?
In 1981, Erd\"{o}s conjectured that for any fixed $0<\alpha<\pi$, the chromatic numbers $\chi(G(S^n;\alpha))$ go to infinity as $n\to\infty$~\cite{erdos1981problem}.
This conjecture was proven by Lov{\'a}sz using topological techniques~\cite{lovasz1983self}.
We also refer the reader to Lov{\'a}sz' famous proof of the chromatic numbers of Kneser graphs~\cite{Lovasz1978}, which gave birth to the field of topological combinatorics, and which mention Borsuk graphs as a more geometric analogue of the combinatorially-defined Kneser graphs.
Since $G(S^n;\alpha)\subseteq \bor{S^n}{\alpha}$, it is therefore also true that $\chi(\bor{S^n}{\alpha})\to\infty$ as $n\to\infty$.
Indeed, Lov{\'a}sz shows that $\chi(\bor{S^n}{\alpha})=n+2$ for all $\alpha$ sufficiently close to $\pi$~\cite{lovasz1983self}.
For the remainder of this paper, we will focus attention on $\bor{S^n}{\alpha}$ instead of $G(S^n;\alpha)$.

Raigarodskii~\cite{raigorodskii2010chromatic, raigorodskii2012chromatic} points out that one of Lov{\'a}sz' claims in~\cite{lovasz1983self} is not correct (even though his result $\chi(\bor{S^n}{\alpha})\to\infty$ still is)\footnote{Raigarodskii furthermore shows that $\chi(\bor{S^n}{\alpha})$ grows exponentially fast in $n$, for any fixed $0<\alpha<n$.}.
Indeed, let $r_n=\arccos{(-\frac{1}{n+1})}$ be the diameter of a regular $(n+1)$-dimensional simplex inscribed in $S^n$.
Also, let $s_n$ be the diameter of the \emph{radial projection} of a single $n$-dimensional face of a regular $(n+1)$-dimensional simplex inscribed in $S^n$.
From~\cite{SantaloConvexRegionsSphericalSurface} (after plugging $\cos \ell=-\frac{1}{n+1}$ into (2.17) and (2.18)), we have
\[s_n=\begin{cases} 
\arccos{(-\frac{n+1}{n+3})} & \text{for }n \text{ odd} \\
\arccos{(-\sqrt{\frac{n}{n+4}})} &  \text{for }n \text{ even}.
\end{cases}\]

\begin{figure}[htb]
\begin{center}
\includegraphics[width=6.5in]{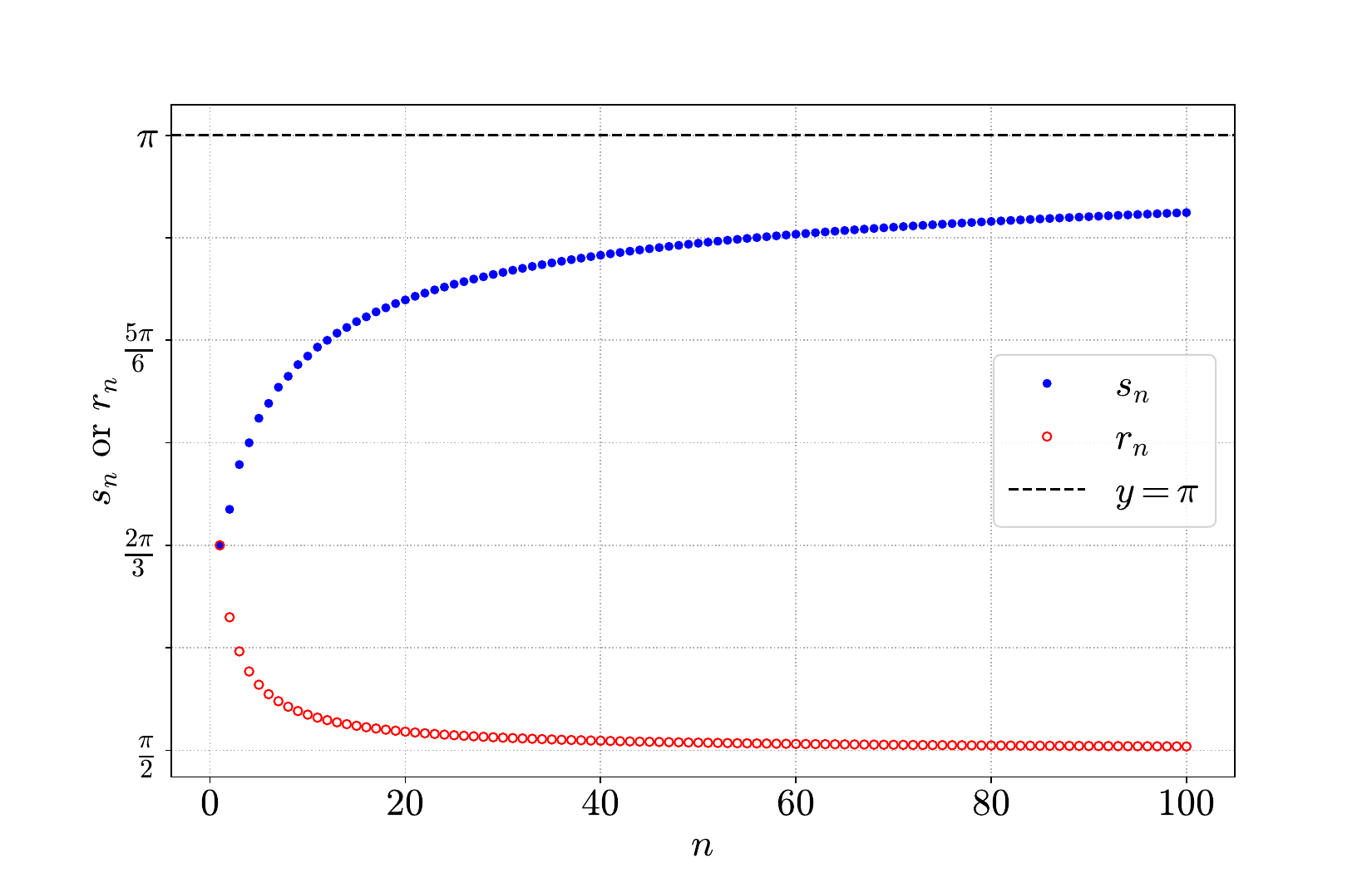}
\caption{Plot of $s_n$ and $r_n$ vs.\ $n$ for $n=1$ to $100$, where $r_n$ is the diameter of the vertex set of an inscribed regular $(n+1)$-simplex in $S^n$, and where $s_n$ is the diameter of a radially projected $n$-face.}
\label{fig:s_n}
\end{center}
\end{figure}

Note $r_1=s_1$ and $r_n<s_n$ for all $n\ge 2$.
Also, $s_n$ is an \emph{increasing} function of $n$ (see Figure~\ref{fig:s_n}), and we have $s_n\to \pi$ as $n\to \infty$. 
Though Lov{\'a}sz claimed $\chi(\bor{S^n}{\alpha})=n+2$ for all $r_n<\alpha<\pi$, Raigorodskii shows that his proof technique in fact only gives $\chi(\bor{S^n}{\alpha})=n+2$ for all $s_n<\alpha<\pi$~\cite{raigorodskii2010chromatic, raigorodskii2012chromatic}.

For any graph $G=(V,E)$, let $N(G)$ be the \emph{neighborhood simplicial complex} of $G$, defined as follows.
The vertex set of $N(G)$ is $V$, i.e.\ the same as the vertex set of $G$.
A finite subset $\{v_0,\ldots,v_k\} \subseteq V$ forms a simplex of $N(G)$ if there is vertex $w \in V$ that is adjacent to each of $v_0, \ldots, v_k$, i.e.\ if $v_0, \ldots, v_k$ share a common neighbor in $V$.
In~\cite[Theorem~1]{Lovasz1978}, Lov{\'a}sz proves that if the first $k$ homotopy groups of $N(G)$ vanish, then the chromatic number of $G$ is at least $k+3$ (see also~\cite[Theorem~5.9.4]{matousek2003using}).
In the case when $G=\bor{S^n}{\alpha}$ is a Borsuk graph, then the neighborhood complex $N(G)$ is a \v{C}ech complex built on $S^n$ at scale $\pi-\alpha$, which for $\frac{\pi}{2}<\alpha<\pi$ is homotopy equivalent to $S^n$ by the nerve lemma; see~\cite{adams2025homotopy}.
Hence for $\frac{\pi}{2}<\alpha<\pi$, the neighborhood complex $N(\bor{S^n}{\alpha})$ has its first $n-1$ homotopy groups vanishing, and the chromatic number of $\bor{S^n}{\alpha}$ is at least $n+2$.

\begin{lemma}[\cite{raigorodskii2010chromatic, raigorodskii2012chromatic,lovasz1983self}]
\label{lem:ChromaticNumberBorsukGraph}
$\chi(\bor{S^n}{\alpha})= n+2$ for $s_n<\alpha<\pi$.
\end{lemma}

\begin{proof}
As per the last paragraph, $\chi(\bor{S^n}{\alpha})\ge n+2$ for $\alpha<\pi$.
We only need to show that $\chi(\bor{S^n}{\alpha})\le n+2$ for $\alpha>s_n$.
Consider the  regular $n+1$-simplex inscribed $S^n$ and project the $n+2$ faces of it onto $S^n$.
The diameter of such image of the faces is $s_n$.
Since $\alpha> s_n$, the vertices of any edge in $\bor{S^n}{\alpha}$ are more than $s_n$ apart.
Thus if we use different colors for vertices in each different such set, we get an $(n+2)$-coloring of $\bor{S^n}{\alpha}$.
\end{proof}

How do the chromatic numbers of $\bor{S^n}{\alpha}$ change as the value of $0<\alpha<s_n$ decreases?
Some asymptotics for the chromatic number of $\bor{S^n}{\alpha}$ as the scale $\alpha$ decreases towards zero are given in~\cite{MoyThesis}, but we are also interested in exact values.

\subsection{Circular Chromatic Numbers}
\label{ssec:circular-chromatic}
Borsuk graphs are naturally related to chromatic numbers of graphs via a related quantity, the circular chromatic number
All of the material in this section on circular chromatic numbers is standard, and can be found in~\cite{CIRCULAR_CHROMATIC_ZHU2001371,HellNesetril2004}.

The circular chromatic number of a graph $G = (V(G), E(G))$ can be defined as follows.
For $r \geq 1$, an \emph{$r$-coloring} of a graph $G$ is a function $\psi\colon V(G)\to S^1$ such that $d(\psi(v),\psi(v')) \geq \frac{2\pi}{r}$ whenever $v,v'$ are adjacent in $G$ (here, $d$ denotes the geodesic distance on $S^1$).
Then the \emph{circular chromatic number of $G$} is
\[ \chi_C(G) \coloneqq \inf\{r \geq 1 \ | \ \textup{there exists an $r$-coloring of $G$}\}.\]

\begin{lemma} \label{lem:circ_chrom_number_borsuk_def} For any graph $G$,
\begin{align*}
    \chi_C(G) & = \inf\{r \geq 1 \ | \ \textup{there exists a graph homomorphism $\phi:G\to \bor{S^1}{\tfrac{2\pi}{r}}$}\}\\
    & = \inf\{\tfrac{2\pi}{\alpha} \geq 1 \ | \ \textup{there exists a graph homomorphism $\phi:G\to \bor{S^1}{\alpha}$}\}.
\end{align*}
\end{lemma}

\begin{proof}
If the graph $G$ is $r$-colorable, then there exists a function $\psi\colon V(G)\to S^1$ such that $d(\psi(v),\psi(v')) \geq \frac{2\pi}{r}$ whenever $v,v'$ are adjacent in $G$.
Thus, $\psi$ is a graph homomorphism from $G$ to $\bor{S^1}{\tfrac{2\pi}{r}}$.
The converse follows similarly.
Hence, we get the first equality.
The second equality is clear after replacing $\frac{2\pi}{r}$ with $\alpha$.
\end{proof}

Recall that, $n$-coloring of a graph $G$ can be viewed as a graph homomorphism $G \to K_n$, where $K_n$ is the complete graph on $n$ vertices.
An $r$-coloring replaces $K_n$ with the Borsuk graph $\bor{S^1}{\frac{2\pi}{r}}$ by Lemma~\ref{lem:circ_chrom_number_borsuk_def}.
The name of the circular chromatic number is justified by its relationship to the usual chromatic number, given in the following lemma.

\begin{lemma}
\label{lemma:chi=ceiling(chi_C)}
For any graph $G$ with finite chromatic number, we have $\chi(G) = \lceil \chi_C(G) \rceil$.
\end{lemma}

\begin{proof}
Let $n \geq 1$.
It is sufficient to show there exists a pair of graph homomorphisms $K_n \to \bor{S^1}{\frac{2\pi}{n}}$ and $\bor{S^1}{\frac{2\pi}{n}} \to K_n$.
Indeed, this will show there exists a graph homomorphism $G \to K_n$ if and only if there exists a graph homomorphism $G \to \bor{S^1}{\frac{2\pi}{n}}$, and thus $\chi_C(G) \leq n$ if and only if $\chi(G) \leq n$.
    
Let $\{ 0, 1, \dots, n-1 \}$ be the vertices of $K_n$ and write points on the circle as angles in $[0, 2 \pi)$.
A homomorphism $K_n \to \bor{S^1}{\frac{2\pi}{n}}$ is given by $k \mapsto \frac{2 \pi k}{n}$, since two distinct points on the circle of the form $\frac{2 \pi k}{n}$ and $\frac{2 \pi k'}{n}$ are at a distance of at least $\frac{2 \pi}{n}$.
A homomorphism $\bor{S^1}{\frac{2\pi}{n}} \to K_n$ is given by sending any point on the circle in $[\frac{2 \pi k}{n}, \frac{2 \pi (k+1)}{n})$ to $k$, since then any two points at distance greater than $\frac{2 \pi}{n}$ are sent to different vertices in $K_n$.
\end{proof}

\section{Nonexistence of graph homomorphisms between Borsuk graphs}
\label{sec:no-graph-homomorphisms}

In this section we will prove a nonexistence result regarding graph homomorphisms between Borsuk graphs.
Recall from Definition~\ref{def:graph-homomorphism} that a graph homomorphism is a function between vertex sets that sends edges to edges.
Therefore, a graph homomorphism $f\colon \bor{S^k}{r} \to \bor{S^n}{\alpha}$ is a function $S^k\to S^n$ that sends points at distance $\ge r$ to points at distance $\ge \alpha$.
We prove in Theorem~\ref{thm:no-graph-homomorphism} that if $f\colon \bor{S^k}{r} \to \bor{S^n}{\alpha}$ is a graph homomorphism with $k>n$ and $r<\pi$, then necessarily $\alpha \le \frac{2\pi}{3}$.
We prove the result using the topology of anti-Vietoris--Rips thickenings.

A weaker version of this result can be deduced from the known chromatic numbers of Borsuk graphs of spheres at scale parameters bigger than $s_n$.
Recall that $\chi(\bor{S^k}{r}) \ge k+2$ for $r<\pi$.
Also, from Lemma~\ref{lem:ChromaticNumberBorsukGraph}, 
$\chi(\bor{S^n}{\alpha}) = n+2$ when $s_n <\alpha <\pi$.
A graph homomorphism $f\colon \bor{S^k}{r} \to \bor{S^n}{\alpha}$ for $\alpha > s_n$ would imply $k+2\le \chi(\bor{S^k}{r}) \le n+2$, which is impossible since $k>n$.
From that, we can deduce that if such a graph homomorphism exists, then necessarily $\alpha \le s_n$.
Since $s_n$ is an increasing function with $s_1=\frac{2\pi}{3}$ and $s_n\to \pi$ as $n\to \infty$, this conclusion from chromatic numbers is strictly weaker than our Theorem~\ref{thm:no-graph-homomorphism} for $n\ge2$.

First, we want to make an observation related to mapping cylinders that will appear in the proof of this theorem.
The \emph{mapping cylinder $M_f$} of a continuous function $f\colon X\to Y$ is the quotient space
\[ M_f=\frac{(X\times I)\coprod Y}{(x,1)\sim f(x)}.\]
The mapping cylinder $M_f$ is always homotopy equivalent to $Y$, simply by contracting $X\times I$ down to $X\times\{1\}$, which as been identified with its image in $Y$.

As an example, consider the mapping cylinder of the 2-to-1 covering space map $f\colon S^n \to \RP^n$ defined by $f(x)=\{x,-x\}$ for all $x\in S^n$.
In this case, $M_f$ is homeomorphic to $\avrm{S^n}{\pi}$, via the homeomorphism $h\colon M_f \to \avrm{S^n}{\pi}$ defined by $h((x,t))=(1-\tfrac{t}{2})\delta_x+\tfrac{t}{2}\delta_{-x}\in \avrm{S^n}{\pi}$ for all $x\in S^n$ and $t\in I$.\footnote{Since $f$ is surjective, we do not need to specify where the homeomorphism maps points in $Y\subseteq M_f$, since that is already determined by where the homeomorphism maps the points in $X\times\{1\}$.
Indeed, $\{x,-x\}\in Y=\RP^n$ is identified in $M_f$ with both of the points $(x,1),(-x,1)\in S^n\times\{1\}$, and hence gets identified under our homeomorphism with $h((x,1))=h((-x,1))=\tfrac{1}{2}\delta_x+\tfrac{1}{2}\delta_{-x}\in \avrm{S^n}{\pi}.$
}
This observation implies $\avrm{S^n}{\pi}\cong M_f\simeq \RP^n$, giving another proof of Lemma~\ref{lem:piCase}.

We now describe a well-known fact about maps between projective spaces, which we will later use as our topological obstruction to the existence of graph homomorphisms between Borsuk graphs.

\begin{lemma}
\label{lem:non-trivial-on-pi_1}   
Let $k>n$, and let $f\colon \RP^k\to \RP^n$ be a continuous map.
Then the induced map $f_*\colon \pi_1(\RP^k)\to \pi_1(\RP^n)$ is trivial.
\end{lemma}

Though this lemma is well-known, we give a proof for completeness.
Our proof uses the ring structure on cohomology, but more elementary proofs can be obtained using covering spaces~\cite{3669362}.
%(\footnote{\url{https://math.stackexchange.com/q/3669362/801498}})

\begin{proof}
Towards a contradiction, suppose the induced map $f_*\colon \pi_1(\RP^k)\to \pi_1(\RP^n)$ is non-trivial.
For $n\geq 2$, we have $\pi_1(\RP^n)\cong \Z/2$.
Since the fundamental group is abelian, it is isomorphic to the integral (first) homology group by the Hurewicz theorem.
Thus we have an isomorphism $f_\ast\colon \Z/2\to \Z/2$.
Using universal coefficient theorem, we have $H^1(\RP^n;\Z/2)\cong \mathrm{Hom}(H_1(\RP^n;\Z),\Z/2).$
Thus we obtain an isomorphism $f^*\colon H^1(\RP^n;\Z/2)\to H^1(\RP^k;\Z/2)$ or $f^*\colon \Z/2\to \Z/2$.
This contradicts Exercise~3(a) of Section 3.2 of~\cite{Hatcher}, which says there is no map $\RP^k\to\RP^n$ inducing a non-trivial map on $H^1$ if $k>n$.
For a sketch of this exercise, let $H^*(\RP^n;\Z/2)\cong (\Z/2)[\alpha]/(\alpha^{n+1})$ where $\alpha$ is the generator of $H^1(\RP^n;\Z/2)$.
Similarly, $H^*(\RP^k;\Z/2)\cong (\Z/2)[\beta]/(\beta^{k+1})$ where $\beta$ is the generator of $H^1(\RP^k;\Z/2)$.
If $f^*\colon H^1(\RP^n;\Z/2)\to H^1(\RP^k;\Z/2)$ is an isomorphism, then $f^*(\alpha)=\beta$.
But this would give the contradiction $0=f^*(0)=f^*(\alpha^{n+1})=(f^*(\alpha))^{n+1}=\beta^{n+1}\neq 0$.

For the edge case $n=1$, at the level of fundamental groups we have the map $f_*\colon \Z/2\to \Z$, which is already trivial.
Thus for $k>n$, any continuous map $f\colon \RP^k\to \RP^n$ induces the trivial map on fundamental groups.
\end{proof}

The proof of Theorem~\ref{thm:no-graph-homomorphism} will rely on the following lemma.
We remind the reader that total anti-VR thickenings are defined in Definition~\ref{def:TAVR}.

\begin{lemma}
\label{lem:phi-map}
Let $(Y,d)$ be a totally bounded metric space.
For $\varepsilon>0$, let $X$ be a finite $\frac{\varepsilon}{2}$-net in $Y$.
Then for any $r>0$, there exists a continuous map $\phi\colon\tavrm{Y}{r+\varepsilon} \to \tavrm{X}{r}$ such that $\supp\{\phi(\mu)\}\subseteq B(\supp(\mu);\frac{\varepsilon}{2}) = \cup_{z\in\supp(\mu)}B(z;\frac{\varepsilon}{2})$ for any measure $\mu\in \tavrm{Y}{r+\varepsilon}$.
\end{lemma}

\begin{proof}
Since $Y$ is a metric space, there exists a partition of unity $\{\rho_x\colon Y\to I\}_{x\in X}$ subordinate to the open cover $\{B(x;\frac{\varepsilon}{2})\}_{x\in X}$ of $Y$.
So, for all $y\in Y$ we have $\sum_{x\in X}\rho_x(y)=1$ with only finitely many $\rho_x(y)$ nonzero, and we have $\rho_x(y)=0$ if $y\notin B(x;\frac{\varepsilon}{2})$.
Define $\phi\colon\tavrm{Y}{r+\varepsilon} \to \tavrm{X}{r}$ in the following manner.
For $\delta_y\in \tavrm{Y}{r+\varepsilon}$ a measure supported on the singleton set $\{y\}$, we define 
\[\phi(\delta_y)=\sum_{x\in X}\rho_x(y)\delta_x \in \tavrm{X}{r}.\]
Clearly, $\supp(\phi(\delta_y))\subseteq B(y;\frac{\varepsilon}{2})$.
We then extend $\phi$ linearly to all of $\tavrm{Y}{r+\varepsilon}$, meaning 
\[\sum_i\lambda_i\delta_{y_i}+\sum_j\gamma_j\delta_{z_j}\mapsto \sum_i\lambda_i\phi(\delta_{y_i})+\sum_j\gamma_j\phi(\delta_{z_j}).\]
Note that if $\phi(\delta_y)$ has non-zero weight on $x$, then $d(x,y)< \frac{\varepsilon}{2}$.
Thus for a general measure $\mu=\sum_i\lambda_i\delta_{y_i}+\sum_j\gamma_j\delta_{z_j} \in \tavrm{Y}{r+\varepsilon}$, the image $\phi(\mu)$ has support $\supp(\phi(\mu))$ contained in $B(\supp(\mu);\frac{\varepsilon}{2}) \coloneqq \cup_{z\in\supp(\mu)}B(z;\frac{\varepsilon}{2})$.
Since the $\{y_i\}$ and $\{z_j\}$ clusters are at least $r+\varepsilon$ apart, the clusters in the image will be at least $r$ apart.
Hence $\phi$ is well-defined.

We next show that $\phi$ is continuous.
Since $Y$ is totally bounded and hence bounded, by~\cite[Corollary~A.2]{AMMW} or~\cite{bogachev2018weak,gibbs2002choosing}, it suffices to show that for every bounded, continuous function $f\colon Y\to \R$ and for any sequence $(\mu_n)$ converging to $\mu\in \tavrm{Y}{r + \varepsilon}$, we have $\int_Y fd\phi(\mu_n)\to \int_Y fd\phi(\mu)$.
Given a bounded, continuous function $f\colon Y\to \R$, define $\tilde{f}\colon Y\to \R$ by $\tilde{f}(y)\coloneqq \sum_{x\in X}\rho_x(y)f(x)$.
Then $\tilde{f}$ is bounded and continuous and $\int_Y fd\phi(\nu) = \int_Y \tilde{f}d \nu$ for every $\nu\in \tavrm{Y}{r+\varepsilon}$.
Since $\mu_n\to \mu$, we have $\int_Y fd\phi(\mu_n) = \int_Y \tilde{f}d\mu_n \to \int_Y \tilde{f}d\mu = \int_Y fd\phi(\mu)$, as desired.

We recall that a graph homomorphism $f\colon \bor{S^k}{r} \to \bor{S^n}{\alpha}$ is equivalent to a possibly discontinuous function $f\colon S^k \to S^n$ such that $d(x,x')\ge r$ implies $d(f(x),f(x'))\ge\alpha$ for all $x,x'\in S^k$.
\end{proof}

\begin{theorem}
\label{thm:no-graph-homomorphism}
For all $k>n$ and $r<\pi$, there is no graph homomorphism $f\colon \bor{S^k}{r}\to \bor{S^n}{\alpha}$ when $\alpha > \frac{2\pi}{3}$.
\end{theorem}

\begin{proof}
If we can prove Theorem~\ref{thm:no-graph-homomorphism} in the special case $\frac{2\pi}{3}< r < \pi$, then the result follows for all $r < \pi$.
Indeed, for $r \le \frac{2\pi}{3}$, this follows simply by considering the composition $\bor{S^k}{\frac{2\pi}{3}}\hookrightarrow \bor{S^k}{r} \to \bor{S^n}{\alpha}$.
Hence we will prove the theorem for $\frac{2\pi}{3}< r < \pi$.

Suppose there exists a graph homomorphism $f\colon \bor{S^k}{r}\to \bor{S^n}{\alpha}$ with $\frac{2\pi}{3}< r < \pi$ and $\alpha > \frac{2\pi}{3}$; we will derive a topological contradiction.
For expository purposes, for the moment assume that $f \colon S^k\to S^n$ is a \emph{continuous} map.
After explaining this simpler case of the proof, we will then relax this assumption and allow $f\colon S^k\to S^n$ to be an arbitrary function.

If $f$ is continuous, then we get the following diagram of continuous maps:
\[ \RP^k \simeq \avrm{S^k}{r} \xrightarrow{f} \avrm{S^n}{\alpha} \simeq \RP^n.\]
The first homotopy equivalence above is since $\frac{2\pi}{3} < r < \pi$, and the second is since $\frac{2\pi}{3} < \alpha < \pi$; see Theorem~\ref{thm:avrmSn-homotopy-type}.
The induced map in the middle is continuous by Lemma~\ref{lem:induced-continuous}.

Let $x\in S^k$ and let
$\omega\colon I\to S^k$ be a geodesic in $S^k$ from $-x$ to $x$.
Consider the following loop $\gamma$ in $\avrm{S^k}{r}$ starting and ending at $\delta_x$, depicted in blue in $\avrm{S^k}{r}$
in Figure~\ref{fig:f_cts_loop}(top left).
\begin{equation}
\label{eqn:gamma_loop}
\gamma(t)=\begin{cases} 
      (1-2t)\delta_x+2t\delta_{-x} & 0\leq t\leq \frac{1}{2} \\
       \delta_{\omega(2t-1)}& \frac{1}{2}\leq t\leq 1,
   \end{cases}    
\end{equation}
Now, $f$ maps $\gamma$ to a loop $f(\gamma)$ in $\avrm{S^n}{\alpha}$ going through $\delta_{f(x)}$ and $\delta_{f(-x)}$.
\[f(\gamma(t))=\begin{cases} 
      (1-2t)\delta_{f(x)}+2t\delta_{f(-x)} & 0\leq t\leq \frac{1}{2} \\
       \delta_{f(\omega(2t-1))}& \frac{1}{2}\leq t\leq 1.
   \end{cases}
\]
Note that $d(f(x),f(-x))\ge\alpha$, i.e., $f(x)$ and $f(-x)$ are nearly antipodal in $S^n$.

We will use the flashlight homotopy as in the proof of Theorem~\ref{thm:avrmSn-homotopy-type},
namely $H\colon \avrm{S^n}{\alpha}\times I\to\avrm{S^n}{\alpha}$ with $h_0\coloneqq H(\,\cdot\,,0)=id_{\avrm{S^n}{\alpha}}$ and $h_1\coloneqq H(\,\cdot\,,1)\colon \avrm{S^n}{\alpha}\to \avrm{S^n}{\pi}$, which exists since $\alpha>\frac{2\pi}{3}$.
Note $h_1(f(\gamma))$ is a loop in $\avrm{S^n}{\pi}$,
\[h_1(f(\gamma(t)))=\begin{cases} 
\delta_{\omega_{f(x),y}(t)} & 0\leq t\leq t_* \\
\frac{1}{1/2-2t_*}[(t-t_*)\delta_{-y}+(1/2-t_*-t)\delta_{y}] & t_*\leq t\leq \frac{1}{2}-t_* \\
\delta_{\omega_{-y,f(-x)}(t)} & \frac{1}{2}-t_*\leq t\leq \frac{1}{2} \\
\delta_{f(\omega(2t-1))}& \frac{1}{2}\leq t\leq 1,
\end{cases}
\]
where we need to explain this notation.
Here $y$ and $-y$ are the antipodal pair whose diameter is parallel to the line joining $f(x)$ and $f(-x)$ (in Figure~\ref{fig:flashlight_homotopy}, compare $f(x),f(-x), y, -y$ with $x_0,x_0',x_1,-x_1$, respectively).
Here $\omega_{f(x),y}\colon [0,t_*]\to S^n$ is the geodesic from $f(x)$ to $y$, and $\omega_{-y,f(-x)}\colon [\frac{1}{2}-t_*, \frac{1}{2}]\to S^n$ is the geodesic from $-y$ to $f(-x)$ on $S^n$, where we choose $0<t_*<\frac{1}{4}$.

We then use the projection to the central core map in Lemma~\ref{lem:piCase}, $(1-t)\delta_x+t\delta_{-x}\mapsto \frac{1}{2}(\delta_x+\delta_{-x})$, to get a loop in the central core (which is identified with $\RP^n$, i.e., we identify the point $\frac{1}{2}(\delta_x+\delta_{-x})$ with $\{\pm x\}\in\RP^n$).
We call this loop $\gamma^\prime$.
\[\gamma'(t)=\begin{cases} 
      \{\pm\omega_{f(x),y}(t)\} & 0\leq t\leq t_* \\
      \{\pm y\} & t_*\leq t\leq \frac{1}{2}-t_* \\
      \{\pm\omega_{-y,f(-x)}(t)\} & \frac{1}{2}-t_*\leq t\leq \frac{1}{2} \\
       \{\pm f(\omega(2t-1))\}& \frac{1}{2}\leq t\leq 1.
   \end{cases}
\]

Choose $\{\pm y\}\in \RP^n$ on $\gamma^\prime$ to be the base point.
Reparametrize $\gamma'\colon [0,1]\to \RP^n$ as $\gamma'\colon [t_*,1+t_*]\to \RP^n$ by taking all inputs modulo $1$, so that $\gamma'$ now starts and ends at $\{\pm y\}$, shown as the blue loop in Figure~\ref{fig:f_cts_loop}(bottom left).
We claim that when we lift the loop $\gamma^\prime$ to the covering space $S^n$ of $\RP^n$, we get a \emph{path} from $-y$ to $y$, or vice versa.
(Recall that lifts of paths are unique once a lift of the basepoint has been chosen.)
For instance, one of the lifts $\tilde{\gamma}\colon [t_*, 1+t_*]\to S^n$ has the form
\[\tilde{\gamma}(t)=\begin{cases} 
      -y & t_*\leq t\leq \frac{1}{2}-t_* \\
      \omega_{-y,f(-x)}(t) & \frac{1}{2}-t_*\leq t\leq \frac{1}{2} \\
      f(\omega(2t-1)) & \frac{1}{2}\leq t\leq 1\\
      \omega_{f(x),y}(t) & t_*\leq t\leq 1+t_*,
   \end{cases}
\]
which is shown as the green path in Figure~\ref{fig:f_cts_loop}(bottom left).
The other lift (in red) is a path from $y$ to $-y$ in $S^n$.
This proves that the loop $\gamma'$ based at $\{\pm y\}\in \RP^n$ is non-trivial in $\pi_1(\RP^k)$ by~\cite[Lemma~54.6(c)]{munkres}, since $S^n$ is the universal cover of $\RP^n$.

Since the loop $\gamma'$ in $\RP^n$ lifts to a path between antipodal points in $S^n$, this means that $f$ induces a continuous map from $\RP^k\to \RP^n$ for $k>n$ that is non-trivial on $\pi_1$.
This contradicts the algebraic fact that any continuous map $\RP^k\to \RP^n$ with $k>n$ must be trivial on $\pi_1$ (Lemma~\ref{lem:non-trivial-on-pi_1}).
Hence $\alpha$ can be at most $\frac{2\pi}{3}$.

\begin{figure}[h]
\def\svgwidth{6.5in}
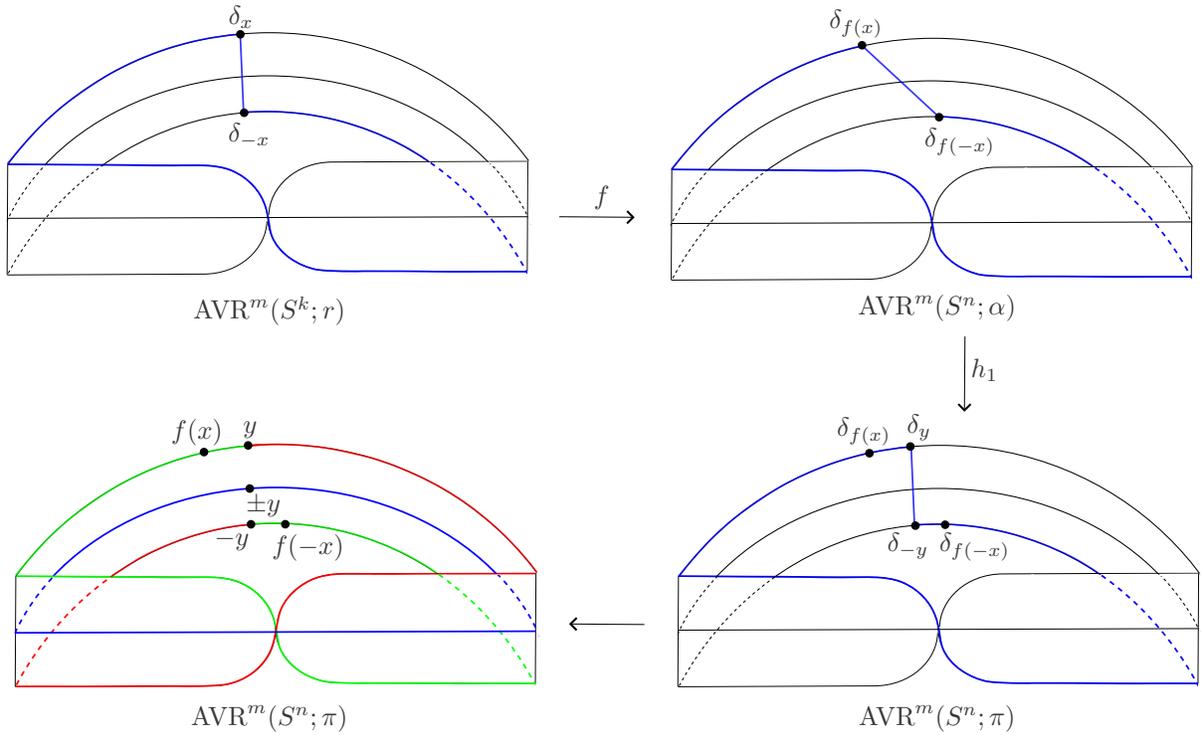
\caption{Non-trivial loop when $f$ is continuous.
The last arrow is the projection to the central core in Lemma~\ref{lem:piCase}.}
\label{fig:f_cts_loop}
\end{figure}

\begin{figure}[htb]
\def\svgwidth{6.5in}
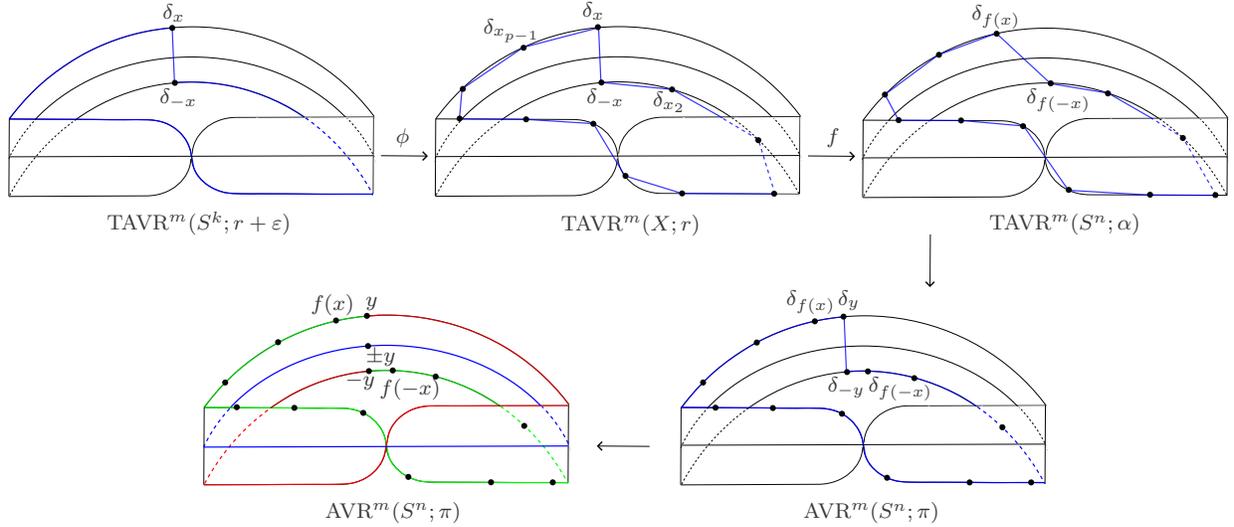
\caption{Non-trivial loop when $f$ may not be continuous.
The second-to-last arrow is the map $\rho$ as in Proposition~\ref{prop:TAVR-to-AVR}, followed by the flashlight homotopy $h_1$.
The last arrow is the projection to the central core.
}
\label{fig:f_discts_loop}
\end{figure}

Now, let us consider the general case when $f$ need not be continuous.
Suppose for a contradiction that there exists a graph homomorphism $f\colon \bor{S^k}{r}\to \bor{S^n}{\alpha}$ with $\frac{2\pi}{3}< r < \pi$ and $\alpha > \frac{2\pi}{3}$.
Choose $\varepsilon>0$ such that $r+\varepsilon<\pi$ and choose a finite $X\subseteq S^k$ that is an $\frac{\varepsilon}{2}$-net, i.e., any point on $S^k$ is within distance at most $\frac{\varepsilon}{2}$ from some point in $X$.
Then we have the following diagram of continuous maps:
\[\RP^k \simeq \tavrm{S^k}{r+\varepsilon} \xrightarrow{\phi} \tavrm{X}{r}\xrightarrow{f}\tavrm{S^n}{\alpha}\simeq \RP^n.\]
The first homotopy equivalence above is since $\frac{2\pi}{3} < r < \pi$, and the last one is since $\frac{2\pi}{3} < \alpha < \pi$ (see Proposition~\ref{prop:TAVR-to-AVR}).
The first map $\phi$ exists by the partition of unity argument in Lemma~\ref{lem:phi-map}.
The second map is induced by $f$, as follows:
$\sum_i\lambda_i\delta_{x_i}+\sum_j\lambda_j\delta_{y_j}\in \tavrm{X}{r}$, where $d(x_i,y_j)\ge r\ \forall i,j$, is sent to $\sum_i\lambda_i\delta_{f(x_i)}+\sum_j\lambda_j\delta_{f(y_j)}\in \tavrm{S^n}{\alpha}$.
As $X$ is finite, $f$ induces a continuous map from $\tavrm{X}{r}$ to $\tavrm{S^n}{\alpha}$ by Corollary~\ref{cor:continuous-on-finite-subset}.

Now we show that $f\phi$ induces a non-trivial map on $\pi_1$.
Without loss of generality, assume that there exists an antipodal pair $\{x,-x\}\in X$; this can be achieved simply by replacing $X$ with $X\cup(-X)$.
Note that the anti-VR thickening is included in the total anti-VR thickening (see Lemma~\ref{lem:atvr-inclusion}), and 
consider the loop $\gamma$ in $\avrm{S^k}{r+\varepsilon}\subseteq \tavrm{S^k}{r+\varepsilon}$ defined as in Equation~\eqref{eqn:gamma_loop}.
By adjusting the sample $X$ if necessary, we can assume there is a sequence of points $x_1, \dots, x_p\in X$ on $\omega$ with $x_1=-x$ and $x_p=x$ such that for any point $\omega(t)$ on the geodesic $\omega$, the intersection $B(\omega(t);\frac{\varepsilon}{2})\cap X$ is equal to either $\{x_i\}$ for some $1\leq i\leq p$ or to $\{x_i, x_{i+1}\}$ for some $1\leq i \leq p-1$.
We further assume that $-x_i\in X$ for all $1\le i \le p$.
Let $\frac{1}{2}=t_1<t_2<\dots <t_p=1$ be such that $\gamma(t_i)=\delta_{x_i}$.
Then a reparametrized loop $(\phi\circ\gamma)(t)$ in $\tavrm{X}{r}$ has the following form:
\begin{equation*} 
(\phi\circ\gamma)(t)=\begin{cases} 
      (1-2t)\delta_x+2t\delta_{-x} & 0\leq t\leq \frac{1}{2} \\
       \tfrac{t_{i+1}-t}{t_{i+1}-t_i}\delta_{x_i}+\tfrac{t-t_i}{t_{i+1}-t_i}\delta_{x_{i+1}}& t_i\leq t\leq t_{i+1} \text{ for } 1\leq i \leq p-1.
   \end{cases}    
\end{equation*}

After applying $f$, we get that the loop $f\circ\phi\circ\gamma$ in $\tavrm{S^n}{\alpha}$ has the following form:
\begin{equation*} 
(f\circ\phi\circ\gamma)(t)=\begin{cases} 
      (1-2t)\delta_{f(x)}+2t\delta_{f(-x)} & 0\leq t\leq \frac{1}{2} \\
       \tfrac{t_{i+1}-t}{t_{i+1}-t_i}\delta_{f(x_i)}+\tfrac{t-t_i}{t_{i+1}-t_i}\delta_{f(x_{i+1})}& t_i\leq t\leq t_{i+1} \text{ for } 1\leq i \leq p-1.
   \end{cases}    
\end{equation*}

Note that $d(x_i, x_{i+1})<\varepsilon< \pi-r$ and $d(-x_i, x_i)=\pi$.
Thus $d(-x_i,x_{i+1})>r$.
Hence
\begin{equation}
\label{eqn:1}
d(f(-x_i), f(x_{i+1}))\ge\alpha,    
\end{equation}
since $f$ is a graph homomorphism from $\bor{S^k}{r}$ to $\bor{S^n}{\alpha}$.
Also, $d(-x_i, x_i)=\pi > r$, thus 
\begin{equation}
\label{eqn:2}
d(f(-x_i), f(x_i))\ge \alpha.
\end{equation}
From~\eqref{eqn:1} and \eqref{eqn:2}, we conclude there exists a measure in $\tavrm{S^n}{\alpha}$ with non-zero weight on all three of $f(x_i)$, $f(x_{i+1})$, $f(-x_i)$, with $f(x_i)$ and $f(x_{i+1})$ in one cluster and $f(-x_i)$ in the other cluster.

Let $\omega_i\colon [t_i, t_{i+1}]\to S^n$ be the geodesic from $f(x_i)$ to $f(x_{i+1})$ for $1\leq i\leq p-1$.
The map $\rho\colon \tavrm{S^n}{\alpha}\to\avrm{S^n}{\alpha}$ in Proposition~\ref{prop:TAVR-to-AVR} maps the line joining $\delta_{f(x_i)}$ and $\delta_{f(x_{i+1})}$, (i.e., measures of the form $c\delta_{f(x_i)}+(1-c)\delta_{f(x_{i+1})}$) to Dirac delta measures $\delta_{\omega_i(t)}$ on the geodesic $\omega_i$.

Thus after applying $\rho$, we get a loop in $\avrm{S^n}{\alpha}$ as follows:
\begin{equation*} 
(\rho\circ f\circ\phi\circ\gamma)(t)=\begin{cases} 
      (1-2t)\delta_{f(x)}+2t\delta_{f(-x)} & 0\leq t\leq \frac{1}{2} \\
       \delta_{\omega_i(t)}& t_i\leq t\leq t_{i+1} \text{ for } 1\leq i \leq p-1.
   \end{cases}    
\end{equation*}
Once again, after applying the flashlight homotopy $h_1$, we get the following loop in $\avrm{S^n}{\pi}$: 

\begin{equation*}
(h_1\circ\rho\circ f\circ\phi\circ\gamma)(t) = \begin{cases}
\delta_{\omega_{f(x),y}(t)} & 0\leq t\leq t_* \\
\frac{1}{1/2-2t_*}[(t-t_*)\delta_{-y}+(1/2-t_*-t)\delta_{y}] & t_*\leq t\leq \frac{1}{2}-t_* \\
\delta_{\omega_{-y,f(-x)}(t)} & \frac{1}{2}-t_*\leq t\leq \frac{1}{2} \\
\delta_{\omega_i(t)}& t_i\leq t\leq t_{i+1} \text{ for } 1\leq i \leq p-1.
\end{cases}
\end{equation*}

Then we apply the projection to the central core  map to get a loop $\gamma'$ in $\RP^n$ as follows:

\[\gamma'(t)=\begin{cases} 
      \{\pm\omega_{f(x),y}(t)\} & 0\leq t\leq t_* \\
      \{\pm y\} & t_*\leq t\leq \frac{1}{2}-t_* \\
      \{\pm\omega_{-y,f(-x)}(t)\} & \frac{1}{2}-t_*\leq t\leq \frac{1}{2} \\
       \{\pm \omega_i(t)\}& t_i\leq t\leq t_{i+1} \text{ for } 1\leq i \leq p-1.
   \end{cases}
\]
As before, $\gamma'$ can be reprametrized $\gamma' \colon [t_*, 1+t_*] \to \RP^n$ so that it starts and ends at $\{\pm y\}$.
When we lift $\gamma'$ to $S^n$ we get the \emph{path} $\tilde{\gamma} \colon [t_*, 1+t_*] \to S^n$ which has the following form:
\[\tilde{\gamma}(t)=\begin{cases} 
      -y & t_*\leq t\leq \frac{1}{2}-t_* \\
      \omega_{-y,f(-x)}(t) & \frac{1}{2}-t_*\leq t\leq \frac{1}{2} \\
     \omega_i(t) & t_i\leq t\leq t_{i+1} \text{ for } 1\leq i \leq p-1\\
      \omega_{f(x),y}(t) & t_*\leq t\leq 1+t_*.
   \end{cases}
\]
This lift is a path in $S^n$ from $-y$ to $y$, depicted as the green path in Figure~\ref{fig:f_discts_loop}(bottom left).
(The red path is the other lift of $\tilde{\gamma}$.)
This once again proves that $f$ induces a continuous map $\RP^k \to \RP^n$, which is non-trivial on $\pi_1$, contradicting Lemma~\ref{lem:non-trivial-on-pi_1}.
Hence for $k>n$ and $r<\pi$, there is no graph homomorphism $f\colon \bor{S^k}{r}\to \bor{S^n}{\alpha}$ when $\alpha > \frac{2\pi}{3}$.
\end{proof}

\section{Conclusion and open questions}
\label{sec:conclusion}

In this paper, we introduce and study anti-Vietoris--Rips metric thickenings, which are related to the previously studied anti-Vietoris--Rips simplicial complexes~\cite{engstrom2009complexes,Jefferson-AATRNtalk2021,JeffersonAntihomology,JeromeRoehm}.
However, since anti-Vietoris--Rips complexes of bounded manifolds are finite-dimensional, whereas Vietoris--Rips complexes of bounded manifolds are infinite dimensional, the relationship between anti-Vietoris--Rips complexes and metric thickenings has a different flavor than the relationship between Vietoris--Rips complexes and metric thickenings.
We have identified the first homotopy types of the anti-Vietoris--Rips metric thickenings of the $n$-sphere for all $n$, and we bound the covering dimension of anti-Vietoris--Rips metric thickenings of manifolds.

We end with a list of open questions.

\begin{question}
Can we identify a set $C \subsetneq [0,\pi]$ of potential \emph{critical} scale parameters satisfying the ``first Morse lemma'' property that if $(s,s')\cap C=\emptyset$ for $s<s'$, then the inclusion $\avrm{S^n}{r'}\hookrightarrow \avrm{S^n}{r}$ is a homotopy equivalence for all $s<r\le r'<s'$?
\end{question}

\begin{question}
Define $r_k=\arccos(-\frac{1}{k+1})$ to be the diameter of the inscribed regular $(k+1)$-simplex in $S^k$.
Recall that $r_1=\frac{2\pi}{3}$.
Does the homotopy type of $\avrm{S^n}{r}$ change as $r$ passes through the scales $\pi > r_1 > r_2 > \ldots > r_{n-1} > r_n$?
\end{question}

\begin{question}
What is $\conn(\avrm{S^n}{r})$ for $0<r<\frac{2\pi}{3}$?
Here $\conn(Y)$ is the largest dimension $k$ such that the homotopy groups $\pi_i(Y)$ are trivial for $i\le k$.
Is $\conn(\avrm{S^n}{r})$ a non-increasing function of $r>0$?

See Theorem~3.13 and Proposition~4.3 of~\cite{engstrom2009complexes} for a related result.
\end{question}

\begin{question}
Can we show that $\avrm{S^n}{r}$ is not contractible for $r>0$?

One way to do this would be by giving an upper bound on $\conn(\avrm{S^n}{r})$.
Indeed, if $\conn(Y)\le k$, then $\pi_i(Y)$ is nontrivial for some $i\ge k+1$, and hence $Y$ is not contractible.
We note that the equivariant techniques of~\cite[Theorem~3]{ABF2} cannot immediately be applied, since the natural action of $\Z/2$ on $\avrm{S^n}{r}$ is not free for $r\le \pi$ (due to fixed points of the form $\frac{1}{2}\delta_x+\frac{1}{2}\delta_{-x}$.)
\end{question}

\begin{question}
What is $\hdim(\avrm{S^n}{r})$ for $0<r<\frac{2\pi}{3}$?
Here, for a topological space $Y$, $\hdim(Y)$ is the largest dimension $k$ in which the homology group $H_k(Y)$ is nonzero.
Is $\hdim(\avrm{S^n}{r})$ a non-increasing function of $r>0$?
%Idea: vertex stars in barycentric subdivision?
\end{question}

\begin{question}
Is $\avrm{S^n}{\frac{2\pi}{3}}$ homotopy equivalent to a $(2n+1)$-dimensional CW complex?
See Theorem~\ref{thm:hom_type_avrm(S^n;2pi/3)}.
\end{question}

\begin{question}
Assuming so, can we use cellular homology to determine the homology groups of $\avrm{S^n}{\frac{2\pi}{3}}$?
\end{question}

\begin{question}
For $M$ a compact manifold, when is $\avrm{M}{r}$ a CW complex?
\end{question}

\begin{question}
Given $k>n$ and a fixed $r<\pi$, what is the largest value of $\alpha$ such that a graph homomorphism $\bor{S^k}{r}\to \bor{S^n}{\alpha}$ exists?

In Theorem~\ref{thm:no-graph-homomorphism}, we showed that for all $k>n$ and $r<\pi$, there is no graph homomorphism $ \bor{S^k}{r}\to \bor{S^n}{\alpha}$ when $\alpha > \frac{2\pi}{3}$.
But this bound on $\alpha$ does not depend on $k$ or $n$, and it does not decrease as $k-n$ grows, as we would expect it to.
\end{question}

\begin{question}
What are the chromatic numbers of the Borsuk graphs $\bor{S^n}{\alpha}$ and $G(S^n;\alpha)$ the scale parameter $\alpha$ decreases?
We remind the reader that for $\alpha$ sufficiently close to $\pi$, the chromatic number of $\bor{S^n}{\alpha}$ is equal to $n+2$~\cite{Lovasz1978,lovasz1983self,raigorodskii2010chromatic, raigorodskii2012chromatic, lovasz1983self, SantaloConvexRegionsSphericalSurface}.
Some asymptotics for the chromatic number of $\bor{S^n}{\alpha}$ as the scale $\alpha$ decreases towards zero are given in~\cite{MoyThesis}, but we are also interested in exact values.
\end{question}

\begin{question}
\label{ques:n-spherical-chromatic-number}
In light of the relationship between the circular chromatic number and the usual chromatic number (Lemma~\ref{lemma:chi=ceiling(chi_C)}), is there a useful definition of an $n$-th \emph{spherical} chromatic number based graph homomorphisms into $\bor{S^n}{\alpha}$?
Relating such a spherical chromatic number to the usual chromatic number will require a better understanding of the chromatic numbers of $\bor{S^n}{\alpha}$ as $\alpha$ varies, as in the previous question.
\end{question}

\section*{Acknowledgements}

We would like to thank Chris Peterson for interesting conversations related to Question~\ref{ques:n-spherical-chromatic-number}.

%\bibliographystyle{plain}
%\bibliography{antiVR.bib}

\appendix

\section{Proof of Lemma~\ref{lem:ExtnOfOptimalTransport}}
\label{app:leftover-proofs}

We give the proof of Lemma~\ref{lem:ExtnOfOptimalTransport}.

\begin{proof}[Proof of Lemma~\ref{lem:ExtnOfOptimalTransport}]
Let $\sigma_1$ and $\sigma_2$ be the marginals of $\sigma$, i.e., $\sigma_1(E)\coloneqq \sigma(E \times X)$ and $\sigma_2(E)\coloneqq \sigma(X \times E)$ for every Borel set $E\subseteq X$.
Thus
\[\mass(\sigma)=\int_{X \times X} d\sigma= \sigma(X \times X)= \sigma_1(X)= \sigma_2(X).\]
Note that for any Borel set $E\subseteq X$, we have $\sigma_1(E)\leq \mu(E)$ and $\sigma_2(E)\leq \nu(E)$.
Hence $\mu-\sigma_1$ and $\nu-\sigma_2$ are Radon measures on $X$.
Define a (full) transport plan between $\mu$ and $\nu$ via $\tilde{\sigma} = \sigma +\tfrac{1}{1-\mass(\sigma)}(\mu-\sigma_1)\times(\nu-\sigma_2)\in \mathcal{P}(X \times X).$
For any Borel set $\tilde{E}\subseteq X\times X$ we have $\sigma(\tilde{E})\leq \tilde{\sigma}(\tilde{E})$, and thus $\tilde{\sigma}$ is an extension of $\sigma$.
Furthermore
\begin{align*}
\mass(\tilde{\sigma})&=\mass(\sigma)+\tfrac{1}{1-\mass(\sigma)}\int_{X\times X}d((\mu-\sigma_1)\times(\nu-\sigma_2))\\
&=\mass(\sigma)+\tfrac{1}{1-\mass(\sigma)}(\mu-\sigma_1)(X)\cdot(\nu-\sigma_2)(X)\\
&= \mass(\sigma)+\tfrac{1}{1-\mass(\sigma)}(\mu(X)-\sigma_1(X))\cdot(\nu(X)-\sigma_2(X))\\
&= \mass(\sigma)+\tfrac{1}{1-\mass(\sigma)}(1-\mass(\sigma))\cdot(1-\mass(\sigma))\\
&= \mass(\sigma)+(1-\mass(\sigma)) = 1.
\end{align*}
Also, for any Borel set $E\subseteq X$, 
\begin{align*}
\tilde{\sigma}(E \times X)&= \sigma(E \times X)+\tfrac{1}{1-\mass(\sigma)}[(\mu-\sigma_1)\times(\nu-\sigma_2)](E \times X)\\
&= \sigma_1(E)+ \tfrac{1}{1-\mass(\sigma)}(\mu-\sigma_1)(E)\cdot (\nu-\sigma_2)(X)\\
&= \sigma_1(E)+ \tfrac{1}{1-\mass(\sigma)}(\mu(E)-\sigma_1(E))\cdot (1-\mass(\sigma))\\
&= \sigma_1(E)+ \mu(E)-\sigma_1(E) = \mu(E).
\end{align*}
Similarly, for any Borel set $E\subseteq X$, we have $\tilde{\sigma}(X \times E)=\nu(E)$.
Hence, $\tilde{\sigma}$ is a transport plan between $\mu$ and $\nu$.
Finally,
\begin{align*}
\cost(\tilde{\sigma})&=\cost\left(\sigma+\tfrac{1}{1-\mass(\sigma)}(\mu-\sigma_1)\times(\nu-\sigma_2)\right)\\
&= \cost(\sigma)+ \tfrac{1}{1-\mass(\sigma)}\cost((\mu-\sigma_1)\times(\nu-\sigma_2))
\end{align*}
with
\begin{align*}
\cost((\mu-\sigma_1)\times(\nu-\sigma_2))&= \int_{X \times X}d(x,y)d((\mu-\sigma_1)\times(\nu-\sigma_2))\\
&\leq \int_{X \times X}\diam(X)d((\mu-\sigma_1)\times(\nu-\sigma_2))\\
&= \diam(X) \int_{X \times X}d((\mu-\sigma_1)\times(\nu-\sigma_2))\\
&= \diam(X) (\mu-\sigma_1)(X)\cdot (\nu-\sigma_2)(X)\\
&= \diam(X) (\mu(X)-\sigma_1(X))\cdot (\nu(X)-\sigma_2(X))\\
&= \diam(X) (1-\mass(\sigma))\cdot (1-\mass(\sigma))\\
&= \diam(X) \left(1-\mass(\sigma)\right)^2.
\end{align*}
Hence,
\[\cost(\tilde{\sigma})\leq\cost(\sigma)+ \tfrac{1}{1-\mass(\sigma)}\diam(X) \left(1-\mass(\sigma)\right)^2= \cost(\sigma)+ \left(1-\mass(\sigma)\right)\diam(X).\qedhere\]
\end{proof}

\section{Generalization of Lemma~\ref{lem:induced-continuous}}
\label{app:gen-induced-continuous}
Here we relax the assumptions in Lemma~\ref{lem:induced-continuous} and give a proof in terms of matchings.

\begin{lemma}
\label{lem:gen-induced-continuous}
Let $X$, $Y$ be metric spaces, and let $K$, $L$ be simplicial complexes with vertex sets $V(K)=X$, $V(L)=Y$.
Let $f\colon X\to Y$ be a map of metric spaces such that the induced map $\tilde{f} \colon K^m \to L^m$ on metric thickenings exists.
If $f$ is continuous and if the image of $f$ has finite diameter, then~$\tilde{f}$ is continuous.
\end{lemma}

\begin{proof}
Since the image of $f$ has bounded diameter, let $C$ be such that $d_Y(f(x),f(x'))\leq C$ for all $x,x'\in X$.
Fix a point $\sum \lambda_i \delta_{x_i}\in K^m$ and fix $\varepsilon>0$.
Using the continuity of $f$ at the finitely many points $x_1,\ldots,x_n$, choose $\nu>0$ so that if $x'\in X$ satisfies $d_X(x_i,x')\le\nu$ for any $i$, then $d_Y(f(x_i),f(x'))\leq \frac{\varepsilon}{2}$.
Reducing $\nu$ if necessary, we can also assume $\nu\le\frac{\varepsilon}{2C}$.
We will show that $W_1(\sum \lambda_i \delta_{x_i},\sum \lambda'_j \delta_{x_j'}) < \nu^2$ implies $W_1(\tilde{f}(\sum \lambda_i \delta_{x_i}),\tilde{f}(\sum \lambda'_j \delta_{x_j'}))\le\varepsilon$, which proves the continuity of $\tilde{f} \colon K^m \to L^m$ at $\sum \lambda_i \delta_{x_i}$.

Let $\pi_{i,j}$ be a matching from $\sum \lambda_i \delta_{x_i}$ to $\sum \lambda'_j \delta{x_j'}$ with $\sum_{i,j} \pi_{i,j}d_X(x_i,x_j')\le\nu^2$.
Let $A=\{(i,j)~|~d_X(x_i,x_j')\geq \nu\}$ and $B=\{(i,j)~|~d_X(x_i,x_j')< \nu\}$.
We have
\[\nu\sum_A\pi_{i,j}\leq \sum_A \pi_{i,j}d(x_i,x_j')\leq \sum_{i,j}\pi_{i,j}d(x_i,x_j')\le\nu^2,\]
so $\sum_A \pi_{i,j}\le\nu$.
Hence
\begin{align*}
W_1\Bigl(\tilde{f}(\textstyle{\sum} \lambda_i \delta_{x_i}),\tilde{f}(\textstyle{\sum} \lambda'_j \delta_{x_j'})\Bigr)
&=W_1\Bigl( \sum_i\lambda_i\delta_{f(x_i)},\sum_j\lambda'_j\delta_{f(x'_j)} \Big)
=W_1\Bigl( \sum_{i,j}\pi_{i,j}\delta_{f(x_i)},\sum_{i,j}\pi_{i,j}\delta_{f(x'_j)} \Bigr)\\
&\le \sum_{i,j} \pi_{i,j}\ d_Y(f(x_i),f(x'_j))\\
&=\sum_A \pi_{i,j}\ d_Y(f(x_i),f(x'_j)) + \sum_B \pi_{i,j}\ d_Y(f(x_i),f(x'_j))\\
&\le C\sum_A \pi_{i,j}+\frac{\varepsilon}{2}\sum_B \pi_{i,j}
\le C\nu+\tfrac{\varepsilon}{2}
\le\varepsilon.
\qedhere
\end{align*}
\end{proof}

\end{document}